\documentclass[a4paper]{article}
\usepackage{fullpage}
\usepackage{amsmath,amsfonts,amsthm,amssymb}
\usepackage{color,xcolor}
\usepackage{ifpdf}
\usepackage{psfrag}
\usepackage{graphicx,graphics,subfigure}
\usepackage{url}
\usepackage{hyperref}
\usepackage{todonotes}

\newtheorem{theorem}{Theorem}[section]

\newtheorem{proposition}[theorem]{Proposition}
\newtheorem{remark}[theorem]{Remark}

\newcommand{\dpar}[2]{\dfrac{\partial #1}{\partial #2}}

 \newcommand{\R}{\mathbb R}

\renewcommand{\P}{\mathbb P}

\newcommand{\bbu}{\mathbf{u}}

\newcommand{\ba}{\mathbf{a}}

\newcommand{\bd}{\mathbf{d}}

\newcommand{\bbf}{\mathbf{f}}

\newcommand{\bm}{\mathbf{m}}
\newcommand{\bn}{\mathbf{n}}

\newcommand{\br}{\mathbf{r}}

\newcommand{\bu}{\mathbf{u}}
\newcommand{\bv}{\mathbf{v}}

\newcommand{\bx}{\mathbf{x}}
\newcommand{\by}{\mathbf{y}}

\newcommand{\bG}{\mathbf{G}}

\newcommand{\bJ}{\mathbf{J}}

\begin{document}
\title{Conservative  scheme compatible with some other conservation laws: conservation of the local angular momentum}
\author{R. Abgrall and Fatemeh Nassajian Mojarrad \\
Institute of Mathematics,
University of Z\"urich\\
Winterthurerstrasse 190, CH 8057 Z\"urich\\
Switzerland}
\date{}
\maketitle
\begin{abstract}
We are interested in building schemes for the compressible Euler equations that are also locally conserving the angular momentum. We present a general framework, describe a few examples of schemes and show results. These schemes can be of arbitrary order.
\end{abstract}
\textbf{keywords.} Euler equations; Angular momentum preservation; Residual distribution schemes; High order methods; Explicit schemes
\tableofcontents

\section{Introduction}
We are interested in the solution of the Euler equations for compressible flows
\begin{equation}
\label{eq:1}
\dpar{\bu}{t}+\text{div }\bbf(\bu)=0,
\end{equation}
defined on a space-time domain $\Omega \times T$, where
$$\bu=(\rho, \rho \bv, E) \text{ and }\bbf=(\rho \bv, \rho \bv\otimes\bv+p\text{Id }, (E+p)\bv)=(\bbf_1, \ldots , \bbf_d), \qquad d=1,2,3.$$
The total energy writes
$$E=\rho e+\frac{1}{2}\rho \bv^2,$$ and the pressure is given by the equation of state
$p=p(\rho, e)$. This function satisfies the usual convexity requirements. Equation \eqref{eq:1} is complemented by initial and boundary conditions.

In addition to \eqref{eq:1}, the solution $\bu$ is known to satisfy other conservation relation. In the case of smooth flows, it is known that $\bu$ satisfies:
\begin{equation}
\label{eq:2}
\dpar{S}{t}+\text{ div }(\rho S)=0 \qquad (\leq 0 \text{ in general}.)
\end{equation}
where $S=-\rho s$ and $s$ is the specific entropy defined by  the Gibbs relation:
$$Tds=d\bigg (\frac{e}{\rho}\bigg )-\frac{p}{\rho^2}d\rho.$$
In order to fully define the entropy, we need a complete equation of state, $p=p(s,\rho)$ from which one can deduce $p=p(\rho,s)$ from standard thermodynamical relations.

There are indeed other relations, for example the preservation of the angular momentum. Defining $\mathbf{J}=\rho \bx\wedge \bv$, and from
$$\dpar{\rho \bv}{t}+\text{ div }(\rho \bv\otimes \bv) +\nabla p=0, $$ we write, in 2D,
$$
\begin{array}{l}
\dpar{\rho \bv}{t}+\dpar{\bbf_1}{x} +\dpar{\bbf_2}{y}=0, \quad \bbf_1=\begin{pmatrix}\rho \bv_x^2+p\\
\rho \bv_x \bv_y\end{pmatrix}=\begin{pmatrix}f_1^1\\f^2_1\end{pmatrix}, \quad \bbf_2=\begin{pmatrix}\rho \bv_x \bv_y\\ \rho \bv_y^2+p\end{pmatrix}=\begin{pmatrix}f_2^1\\f^2_2\end{pmatrix},
\end{array}
$$
We note that $f_1^2=f_2^1$.
\begin{equation*}
\begin{split}
\dpar{\rho\bx\wedge\bv}{t}&=\bx\wedge\dpar{\rho\bv}{t}=-\bx\wedge\big ( \dpar{\bbf_1}{x_1} +\dpar{\bbf_2}{x_2}\big )\\
&=-x_1\big ( \dpar{f_1^2}{x_1}+\dpar{f_2^2}{x_2}\big ) +x_2 \big (  \dpar{f_1^1}{x_1}+\dpar{f_2^1}{x_2}\big )\\
&= -\dpar{}{x_1}\big ( f_1^2 x_1-x_2f_1^1\big )-\dpar{}{x_2}\big (x_1 f_2^2 -f_2^1 x_2\big )\\
&\qquad +\big ( f_1^2-f_2^1\big ),\\
\end{split}
\end{equation*}
Hence we have:
\begin{equation}\label{momentC:2D}
\dpar{\mathbf{J}}{t}+\dpar{}{x_1}\big ( f_1^2 x_1-x_2f_1^1\big )+\dpar{}{x_2}\big (x_1 f_2^2 -f_2^1 x_2\big )=0.
\end{equation}
In 3 dimensions, we have a similar results: a direct computation of the time evolution of the angular momentum equation give:
\begin{equation*}
\begin{split}
\dpar{J_1}{t}&+\bigg ( \dpar{x_2f_1^3}{x_1}+\dpar{x_2f_2^3}{x_2}+\dpar{x_2f_3^3}{x_3}\bigg ) -\bigg ( \dpar{x_3f_1^2}{x_1}+\dpar{x_3f_2^2}{x_2}+\dpar{x_3f_3^2}{x_3}\bigg )\\
& \qquad + \big ( f_3^2-f_2^3\big )=0,\\
\dpar{J_2}{t}&+\bigg ( \dpar{x_3f_1^1}{x_1}+\dpar{x_3f_2^1}{x_2}+\dpar{x_3f_3^1}{x_3}\bigg ) -\bigg ( \dpar{x_1f_1^3}{x_1}+\dpar{x_1f_2^3}{x_2}+\dpar{x_1f_3^3}{x_3}\bigg )\\
& \qquad + \big ( f_1^3-f_3^1\big )=0,\\
\dpar{J_3}{t}&+\bigg ( \dpar{x_1f_1^2}{x_1}+\dpar{x_1f_2^2}{x_2}+\dpar{x_1f_3^2}{x_3}\bigg ) -\bigg ( \dpar{x_2f_1^1}{x_1}+\dpar{x_2f_2^1}{x_2}+\dpar{x_2f_3^1}{x_3}\bigg )\\
& \qquad + \big ( f_2^1-f_1^2\big )=0,
\end{split}
\end{equation*}
and since
$$f_2^3=f_3^2=\rho u_2u_3, \quad  f_3^1=f_1^3=\rho u_3u_1, \quad f_1^2=f_2^1=\rho u_1 u_2,$$
we obtain the additional conservation relation:
\begin{equation}\label{momentC:3D}
\begin{split}
\dpar{J_1}{t}&+ \dpar{}{x_1}\big ( x_2f_1^3- x_3f_1^2\big )+\dpar{}{x_2}\big ( x_2f_2^3-x_3f_2^2\big )+\dpar{}{x_3}\big ( x_2f_3^3-x_3f_3^2\big ) =0,\\
\dpar{J_2}{t}&+\dpar{}{x_1}\big (x_3f_1^1-x_1f_1^3\big ) +\dpar{}{x_2}\big ( x_3f_2^1-x_1f_2^3\big )+\dpar{}{x_3}\big ( x_3f_3^1-x_1f_3^3\big ) =0,\\
\dpar{J_3}{t}&+ \dpar{}{x_1}\big ( x_1f_1^2-x_2f_1^1\big )+\dpar{}{x_2}\big ( x_1f_2^2-x_2f_2^1\big )+\dpar{}{x_3}\big (x_1f_3^2-x_2f_3^1\big ) =0.
\end{split}
\end{equation}
The weak form of the relations \eqref{momentC:2D} and \eqref{momentC:3D} are also satisfied by $\mathbf{J}$. To see this, in the 2D case, one has to take the test functions $x_1\varphi$ and $x_2\varphi$ where $\varphi$ has compact support and is $C^1$ in space-time domain, and subtract  to get the weak form of \eqref{momentC:2D}. The same method also works in 3D.
\begin{remark}\label{translation}
The PDEs \eqref{momentC:2D} and \eqref{momentC:3D} are translational  invariant. This means that if $\mathbf{J}$ satisfies \eqref{momentC:2D} or \eqref{momentC:3D} in the frame of coordinates $\bx$, then $\mathbf{J}'=\mathbf{J}+\rho\mathbf{a}\wedge \bv$ will satisfy \eqref{momentC:2D} or \eqref{momentC:3D} in the frame of coordinates $\bx'=\bx+\mathbf{a}$.
\end{remark}

Developing schemes that locally preserves the angular momentum has recently attracted some attention, see e.g \cite{zbMATH06660499,gaburro} and the reference therein. In \cite{tokareva1,abg,tokareva2,abgrallOeffnerRanocha} has been developed a simple technique that enables to correct an initial scheme so that an additional conservation is satisfied. This has been applied to some non conservative form of the Euler equation, also to satisfy entropy conservation when this is relevant, and to various schemes: residual distribution schemes and discontinuous Galerkin ones. The goal of this paper is to show how one can extend the method to the problem of angular preservation, for second and higher order schemes, in a fully discrete manner, both in time and space.

The paper is organised as follows. In section 2 we briefly introduce a high order residual distribution scheme for the unsteady version of \eqref{eq:1}. In section 3 we construct a second order scheme for \eqref{eq:1} that preserves locally the angular momentum; then we show this is possible for triangles and tetrahedrons. In section 4 we provide a possible solution to conservation of momentum for higher than second order schemes. In section 5 we study the validity and effectiveness of the proposed strategy, considering three different  problems. In section 6 we describe how to extend this to discontinuous Galerkin schemes. Finally, section 7 provide some conclusions and future perspectives.
\section{Numerical scheme}
\subsection{Some generalities}

The goal is to construct a scheme for \eqref{eq:1} that preserves locally the angular momentum. For this we use a high order residual distribution scheme for the unsteady version of \eqref{eq:1} that we briefly describe now, see \cite{abgrall,paola} for more details.

We are given a tessellation of $\Omega$  by non overlapping simplex denoted generically by $K$, and more precisely $\Omega=\cup K$. We assume the mesh to be conformal which means that $K\cap K'$ is either empty, reduce to a full face or a full edge or a vertex.  From this, we consider $\mathcal{V}=\{f:\Omega \rightarrow \R, \quad f_{|K}\in \P^r \quad \forall K \in \Omega\}\cap C^0(\Omega)$. We are given a basis of $\mathcal{V}$, $\{\varphi_\sigma\}$ where the $\sigma$ are the degrees of freedom (DOFs). The restriction of $\varphi_\sigma$ to any $K$ is assumed here to be a B\'ezier polynomial, see appendix \ref{sec:Bezier} for the notations we use. The reason is that for any $\sigma$, we have
\begin{equation}\label{cond:1}|C_\sigma|=\int_{\Omega} \varphi_\sigma\; d\bx>0.\end{equation}
This is not necessarily true for Lagrange polynomials. 

We are looking for an approximation of the solution of the form
$$\bu=\sum\limits_{\sigma}\bu_\sigma \varphi_\sigma,$$
Then the solution $\bu^{n}$ at time $t_n$ can be written as follows
$$\bu^{n}=\sum\limits_{\sigma}\bu^{n}_\sigma \varphi_\sigma.$$
The coefficients $\bu^{n}_\sigma$, are chosen by a numerical method and we use the following technique to calculate the coefficients in the approximation.
\subsection{Residual distribution scheme for steady problems}
First, we consider a steady  version of system \eqref{eq:1}
\begin{equation}
\label{eq:111}
\text{div }\bbf(\bu)=0.
\end{equation}
\begin{figure}[ht]
\begin{center}
\subfigure[Compute the total residual.]{\includegraphics[width=0.3\textwidth]{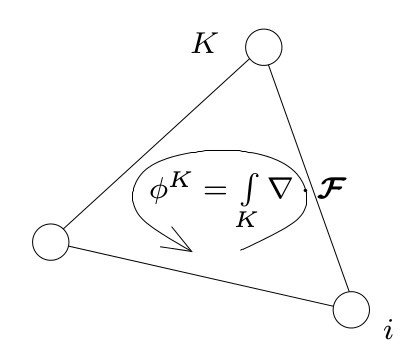}\label{total_residual}}
\subfigure[Compute the nodal residuals.]{\includegraphics[width=0.3\textwidth]{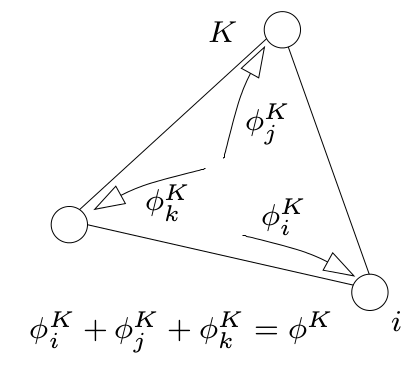}\label{nodal_residuals}}
\subfigure[Collect all the nodal residual contributions.]{\includegraphics[width=0.3\textwidth]{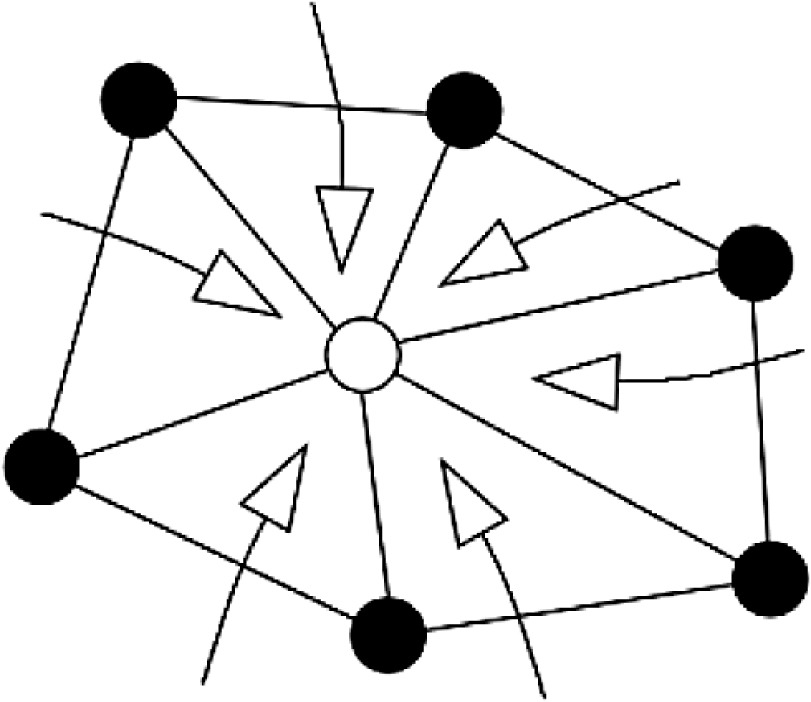}\label{collect}}
\end{center}
\caption{Illustration of the main steps of the residual distribution approach.}
\end{figure}
The main steps of the residual distribution approach can be summarized as follows
\begin{enumerate}
\item For any element $K \in \Omega$, compute a fluctuation term (total residual) (see figure \ref{total_residual})
\begin{equation}
\Phi^K(\bbu)=\int_{\partial K} \bbf(\bu)\cdot \bn\; d\bx \;(=\int_K \text{ div }\bbf(\bu) \; d\bx),
\end{equation}
\item For every DOF $\sigma$ within the element K, define the nodal residuals $\Phi_\sigma ^K$ as the contribution to the fluctuation term $\Phi^K$ (see figure \ref{nodal_residuals}) such that
\begin{equation}
\Phi^K(\bbu)=\sum_{\sigma \in K} \Phi_\sigma ^K, ~\forall K \in \Omega,
\end{equation}
\item The resulting scheme is obtained for each DOF $\sigma$ by collecting
all the nodal residual contributions $\Phi_\sigma ^K$ from all elements K surrounding a node $\sigma \in \Omega$ (see figure \ref{collect})
\begin{equation}
\label{eq:111:2}
\sum_{K,\sigma \in K} \Phi_\sigma ^K(\bbu)=0, ~ \forall \sigma \in \Omega.
\end{equation}
\end{enumerate}
We can write a similar formulation for the boundary conditions. 
For any DOF $\sigma \in \Omega$ we can split \eqref{eq:111} into the internal and boundary contributions
\begin{equation}
\sum_{K,\sigma \in K} \Phi_{\sigma,\bx} ^K(\bbu)+\sum_{\Gamma,\sigma \in \gamma} \Phi_{\sigma,\bx} ^\Gamma(\bbu)=0,
\end{equation}
where $\Phi_{\sigma,\bx} ^K$ and $\Phi_{\sigma,\bx}^ \gamma$ are the residuals
corresponding to the spatial discretization and $\gamma$ is an edge on the boundary $\Gamma$ of $\Omega$. 
If we suppose the Dirichlet boundary condition $\bu=g$ on $\Gamma$, for any K and $\Gamma$, $\Phi_{\sigma,\bx} ^K$ and $\Phi_{\sigma,\bx}^ \gamma$ satisfy the following conservation relations
\begin{subequations}\label{conservation}
\begin{equation}\label{conservation:1}
\sum_{K,\sigma \in K} \Phi_{\sigma,\bx} ^K(\bbu)=\int_{\partial K} \bbf(\bu)\;.\;\bn,
\end{equation}and for the boundary condition (here Dirichlet in weak form)
\begin{equation}\label{conservation:2}
\sum_{\Gamma,\sigma \in \Gamma} \Phi_{\sigma,\bx} ^\Gamma(\bbu)=\int_\Gamma \big( \mathcal{F}_\bn(\bu,g)-\bbf(\bu)\;.\;\bn \big).
\end{equation}
\end{subequations}
In the appendix \ref{appendix:RD}, we provide several examples of such fluctuations. For now, and in order to simplify the discussion on the approximation of the unsteady problem \eqref{eq:1},  we introduce a "variational" form of the fluctuation. In all the known cases, we can write $\Phi_\sigma^K$ as
\begin{equation}\label{eq:beta}
\Phi_\sigma^K(\bbu)=\beta_\sigma^K(\bbu)\; \Phi^K(\bbu),
\end{equation}
where
$$\sum\limits_{\sigma\in K}\beta_\sigma^K(\bbu)=\text{Id} \quad (=1 \text{ in the scalar case}).$$
The distribution coefficients can be scalar in the scalar case, and matrices in the system case.

Then, we can rewrite the scheme in a Petrov-Galerkin fashion:
\begin{equation*}
    \begin{split}
\Phi_\sigma^K(\bbu)&=\beta_\sigma^K(\bbu)\Phi^K(\bbu)=\int_K \beta_\sigma^K(\bbu)\; \text{div }\bbf(\bu)\; d\bx=\int_K\varphi_\sigma \text{ div }\bbf(\bu)\; d\bx+\int_K \big (\beta_\sigma^K(\bu)-\varphi_\sigma \text{Id}\big )\text{ div }\bbf(\bu)\; d\bx\\
&=-\int_K \nabla\varphi_\sigma\cdot\bbf(\bu)\; d\bx+\int_{\partial K}\varphi_\sigma \bbf(\bu)\cdot \bn \;d\gamma+\int_K\xi_\sigma^K(\bbu) \text{ div }\bbf(\bu)\; d\bx,
\end{split}
\end{equation*}
where $\xi_\sigma^K(\bbu)=\beta_\sigma^K(\bbu)-\varphi_\sigma \text{Id}$. To avoid confusion, we will write 
$\xi_\sigma^K(\bbu)$ as $\xi_\sigma^K$, thus removing the \emph{functional} dependency of this term with respect to the solution $\bu$. We remove it, but we do not forget it!

We note that the functions $\xi_\sigma^K$ satisfies
$$\sum\limits_{\sigma\in K}\xi_\sigma^K=0.$$
Then we define $\Psi_{\sigma, \bbu}=\varphi_\sigma \text{Id} +\xi_\sigma$ with, for $\bx\in K$, 
$$\xi_\sigma(\bx)=\xi_\sigma^K.$$
Defining $W=\text{span}(\Psi_\sigma)$, we can formally rewrite the scheme as: for any $v\in \mathcal{V}$, find $\bu$ such that
$$-\int_\Omega\nabla v\cdot \bbf(\bu)\; d\bx+\int_{\partial\Omega} v\mathcal{F}_\bn\; d\gamma+\sum_{K}\int_K \xi(\bbu,v)\; \text{div }\bbf(\bu)\; d\bx=0,$$
where
$$\xi(\bx)=\sum_\sigma v_\sigma \xi_\sigma.$$
Here the functional dependence of $\xi$ with respect to $\bu$ exists but is implicit.
\subsection{Residual distribution scheme for unsteady problems}
Using the results presented above, it becomes possible to describe the unsteady version of the scheme. In order to get a consistent approximation, we simply "multiply" \eqref{eq:1} with the Petrov-Galerkin test function $v+\xi(v)$, and integrate in time. Note that $v$ will not depend on time.
Between $t_n$ and $s\in ]t_n, t_{n+1}]$, we get
$$\int_{t_n}^s \int_\Omega\big (v+\xi(v)\big)\dpar{\bu}{t}\; d\bx\; dt-\int_{t_n}^s\int_\Omega \nabla v\;\cdot \;\bbf(\bu)\; d\bx\; dt
+\int_{t_n}^s\int_{\partial\Omega} 
v\mathcal{F}_\bn\; d\gamma\; dt+\sum_K \int_{t_n}^s\int_K \xi(v) \text{ div }\bbf(\bu)\; d\bx\; dt=0.
$$
Since $v$ is independent of time, this can be equivalently rewritten as:
$$\int_\Omega\big (v+\xi(v)\big) \big ( \bu(\bx,s)-\bu(\bx,t_n)\big ) \; d\bx
-\int_{t_n}^s\int_\Omega \nabla v\;\cdot \;\bbf(\bu)\; d\bx\; dt
+\int_{t_n}^s\int_{\partial\Omega} 
v\mathcal{F}_\bn\; d\gamma\; dt+\sum_K \int_{t_n}^s\int_K \xi(v) \text{ div }\bbf(\bu)\; d\bx\; dt=0,$$
and then, for any $\sigma$,
\begin{equation*}
    \begin{split}
        \sum\limits_{K, \sigma\in K}& \Bigg [ \int_K \big (\varphi_\sigma+\xi_\sigma \big ) \big (\bu(\bx,s)-\bu(\bx,t_n)\big ) \; d\bx +\\
        &\int_{t_n}^s\bigg ( -\int_K\nabla\varphi_\sigma \cdot \bbf(\bu)\; d\bx+\int_K\xi_\sigma\; \text{div }\bbf(\bu)\; d\bx+\int_{\partial K}\varphi_\sigma \bbf(\bbu)\cdot\bn\;d\gamma \bigg )\;dt\Bigg ]+\text{BCs}=0
    \end{split}
\end{equation*}
This suggest to introduce the space-time fluctuations
\begin{equation}
    \label{space-time}
    \Phi_\sigma^K(\bu,s)=\underbrace{\int_K \big (\varphi_\sigma+\xi_\sigma \big ) \big (\bu(\bx,s)-\bu(\bx,t_n)\big ) \; d\bx}_{\Phi_{\sigma,t}^K(\bu,s)}
    +\int_{t_n}^s \Phi_{\sigma,\bx}^K(\bu,s)\;dt
\end{equation}
where
$$\Phi_{\sigma,\bx}^K(\bu,s)=-\int_K\nabla\varphi_\sigma \cdot \bbf(\bu)\; d\bx+\int_K\xi_\sigma\; \text{div }\bbf(\bu)\; d\bx+\int_{\partial K}\varphi_\sigma \bbf(\bbu)\cdot\bn\;d\gamma .$$

\subsection{Iterative timestepping method}
In order to integrate the system, we proceed in two steps. First we introduce sub-time steps in the interval $[t_n, t_{n+1}]$:
 we subdivide the interval $[t_n, t_{n+1}]$ into sub-intervals obtained from a partition
$$t_n=t_{(0)} < t_{(1)} < \ldots < t_{(l)}< \ldots <t_{(M)}=t_{n+1}.$$
In all our examples, this partition will be regular: $t_{(l)}=t_n+\tfrac{l}{M}\Delta t$, but other choices could have been made, especially for very high order accuracy in time.   We construct an approximation of $\bu$ at  times $t_{(l)}$, denoted by 
$\bu_{(l)}\approx \bu(t_{(l)})$. Then we introduce $\mathcal{I}_M$ the Lagrange interpolant (in time) of degree $M$ defined from this subdivision. We also introduce the piece-wise constant interpolant $\mathcal{I}_0$ defined as
$$\mathcal{I}_0(s)= \bu_{(0)} \text{ if } s\in [t_n, t_{n+1}].$$
Another choice could have been made
$$ \mathcal{I}_0(s)= \bu_{(l)} \text{ if } s\in [t_{(l)}, t_{(l+1)}[.$$
 for $l=0, \ldots , M-1$. The notation $U$ represents the vector $U=(\bu_{(0)}, \bu_{(1)}, \ldots , \bu_{(M)})$ i.e the vector of all the approximations for the sub-steps. Note that $\bu_{(0)}=\bu^n$ and $\bu_{(M})=\bu^{n+1}$.
We need residuals, $\{\Phi_{\sigma,t}^{K,p}(U), \Phi_{\sigma,\bx}^{K,p}(U)\}$ computed for time $t_p$ and  that satisfies, for any $p$,
$$\sum_{\sigma\in K}\Phi_{\sigma, t}^{K,p}(U)=\int_K\big (\bu_{(p)}-\bu^n\big )\;  d\bx, \quad \sum_{\sigma\in K}\Phi_{\sigma, \bx}^{K,p}(U)=\int_{t_n}^{t_{(p)} }\int_{\partial K} \mathcal{I}_M(\bbf)\cdot \bn \; d\gamma\;dt.$$
We set $\Phi_\sigma^K(U)=\Phi_{\sigma,t}^{K,p}(U)+\Phi_{\sigma,\bx}^{K,p}(U)$.

Then we introduce two approximations of \eqref{eq:1}:
\begin{itemize}
\item A first order approximation in time: for any $\sigma \in K$ and $l=1, \ldots , M$, 
\begin{equation}\label{cond:2}\big [L_1(U; \bu^n)\big ]_{\sigma, (l)}:=|C_\sigma| \bu_{\sigma, (l)}- |C_\sigma| \bu_{\sigma, (0)}+\sum\limits_{K, \sigma \in K}\int_{t_{(0)}}^{t_{(l)} }
\mathcal{I}_0\big ( \Phi_{\sigma,\bx}^K(U),s\big )\; ds,\end{equation}
where
$$
\mathcal{I}_0\big ( \Phi_{\sigma,\bx}^K(U),s\big )=\mathcal{I}_0\big (\Phi_{\sigma,\bx}^K(\bu_{(0)}),\Phi_{\sigma,\bx}^K(\bu_{(1)}), \ldots ,
\Phi_{\sigma,\bx}^K(\bu_{(M)}),s\big),
$$
i.e. here
\begin{align*}
\big [L_1(U; \bu^n)\big ]_{\sigma, (l)}&=|C_\sigma| \bu_{\sigma, (l)}- |C_\sigma| \bu_{\sigma, (0)}+\big ( t_{(l)}-t_{(0)} \big ) \sum\limits_{K, \sigma \in K}
\mathcal{I}_0\big ( \Phi_{\sigma,\bx}^K(\bu_{(0)})\big )\\
&=|C_\sigma| \bu_{\sigma, (l)}- |C_\sigma| \bu_{\sigma, (0)}+\frac{l}{M}\Delta t  \sum\limits_{K, \sigma \in K}
 \Phi_{\sigma,\bx}^K(\bu_{(0)})
\end{align*}
This is a first order \emph{explicit} approximation in time.
\item A high order approximation: for any $\sigma \in K$ and $l=1, \ldots , M$,  
\begin{equation}\label{cond:3}\big [L_2(U; \bu^n)\big ]_{\sigma, (l)}:=\sum_{K, \sigma\in K} \int_{t_{(0)}}^{t_{(l)}} \mathcal{I}_M\big ( \Phi_\sigma^K (U ),s\big )\; ds,\end{equation}

After having performed exact integration in time
to obtain the approximation for every sub-steps, the time integration can be written in the form
\begin{equation}
\int_{t_{(0)}}^{t_{(l)}} \mathcal{I}_M\big (\Phi_{\sigma,\bx}^K(U)\big)\;ds=\sum_{k=0}^M \theta_{k}^l \; \Phi_{\sigma,\bx}^K (U ).   
\end{equation}
\end{itemize}
Ideally, we would like to solve for each $\sigma$ and each $l$,
\begin{equation}\label{eq:HO}\big[L_2(U; \bu^n)\big ]_{\sigma, (l)}=0,\end{equation}
but this is very difficult in general and the resulting scheme derived by $L_2$ operator is implicit. Instead we use a defect correction method and proceed within the time interval $[t_n, t_{n+1}]$ as follows:
\begin{enumerate}
\item {Set }$U^{(0)}=(\bu^{n},\bu^{n}, \ldots , \bu^{n})$,
\item For any $p\geq 0$, define $U^{(p+1)}$  by
\begin{equation}\label{eq:HO2}
\big [L_1(U^{(p+1)}; \bu^n)\big ]_{\sigma,(l)}=\big [L_1(U^{(p)}; \bu^n)\big ]_{\sigma,(l)}-\big [L_2(U^{(p)}; \bu^n)\big ]_{\sigma,(l)}.\end{equation}
\end{enumerate}
Since $L_1$ is explicit, $U^{(p+1)}$ can be obtained explicitly.
In our case, this amounts to a multi-step method where each step writes as
\begin{equation}
\label{rd:1}
\begin{split}
|C_\sigma| \big (U^{(p+1)}_{\sigma}-U^{(p)}_{\sigma}\big ) &=-  \sum_{K, \sigma\in K} \Phi_\sigma^K(U^{(p)}),
\end{split}
\end{equation}
The scheme \eqref{eq:HO2} is completely explicit. One has the following result, see \cite{abgrall}:
\begin{proposition}
If two operators $L_1$ and $L_2$ depending on the discretization scale $\Delta=\Delta t$, are
such that:
\begin{itemize}
\item There exists a unique $U_{\Delta}^\star$ such that $L_2(U_{\Delta}^\star)=0.$
\item $L_1$ is coercive, i.e., there exists $\alpha_1 > 0$ independent of $\Delta$, such that for any U and V,
$$
\alpha_1 \Vert U-V \Vert\leq \Vert L_1(U)-L_1(V) \Vert,
$$
\item $L_1-L_2$ is uniformly Lipschitz continuous with Lipschitz constant $\alpha_2 \Delta$, i.e., there exists $\alpha_2 > 0$ independent of $\Delta$, such that for any U and V,
$$
\Vert \big(L_1(U)-L_2(U)\big)-\big(L_1(V)-L_2(V)\big) \Vert\leq\alpha_2 \Delta \Vert U-V \Vert.
$$
Then if $\nu = \frac{\alpha_2 \Delta}{\alpha_1} < 1$ the defect correction method is convergent, and after p iterations the error is
smaller than $\nu^p \Vert U^{(0)}-U_{\Delta}^\star\Vert.$
\end{itemize}
\end{proposition}
Therefore if there is a unique solution $U_{\Delta}^\star$ of \eqref{eq:HO}, then after $k$ steps, $\bu^{(k)}_{(M)}=U_{\Delta}^\star+O(\Delta t^k)$, so that only $M+1$ steps are needed. Of course this is true thanks to the conditions \eqref{cond:1}, \eqref{cond:2}, \eqref{cond:3} and the condition that $L_1$ is a first order approximation of $L_2$:
\begin{equation}\label{cond:4}L_1-L_2=O(\Delta t).\end{equation}
\section{Angular momentum preservation: second order case}

It is clear that the residuals  depends where the order in time. 
We can compute the variation of the angular momentum from \eqref{rd:1}. For any $\sigma$, we have
$$
|C_\sigma| \big (\bx_\sigma \wedge(\mathbf{m}^{(p+1)}_\sigma-\mathbf{m}^{(p)}_\sigma)\big ) + \sum_{K, \sigma\in K} \bx_\sigma\wedge\Phi_{\mathbf{m},\sigma}^K(U^{(p)})+\sum_{\Gamma, \sigma\in \Gamma} \bx_\sigma\wedge\Phi_{\mathbf{m},\sigma}^\Gamma(U^{(p)})=0,
$$
so that
$$
|C_\sigma| \big (\mathbf{J}^{(p+1)}_\sigma-\mathbf{J}^{(p)}_\sigma\big ) + \sum_{K, \sigma\in K} \bx_\sigma\wedge\Phi_{\mathbf{m},\sigma}^K(U^{(p)})+\sum_{\Gamma, \sigma\in \Gamma} \bx_\sigma\wedge\Phi_{\mathbf{m},\sigma}^\Gamma(U^{(p)})=0
$$
where here $\Phi_{\mathbf{m},\sigma}^K$ is the momentum component of the residual at the DOF $\sigma$ in the element $K$ and $\Phi_{\mathbf{m},\sigma}^\Gamma$ is the momentum component of the residual at the DOF $\sigma$ in the boundary $\Gamma$.
An easy condition for having local conservation is: 
\begin{equation}
\label{cons:J}
\sum_{\sigma\in K} \bx_\sigma\wedge\Phi_{\mathbf{m},\sigma}^K(U^{(p)})=\int_K \big ( \mathbf{J}^{(p)}-\mathbf{J}^{(0)}\big )d\bx +\Delta t\int_{\partial K} \mathcal{I}_2\big ({\mathbf{G}}(U^{(p)})\big )\cdot \bn\; d\gamma :=\Phi_{\mathbf{J}}^K
\end{equation}
\begin{equation}
\label{cons:JJ}
\sum_{\sigma\in \Gamma} \bx_\sigma\wedge\Phi_{\mathbf{m},\sigma}^\Gamma(U^{(p)})=\Delta t\int_{\partial K} \Big ( \hat{\bG}-\mathcal{I}_2\big ({\mathbf{G}}(U^{(p)})\big)\cdot \bn \Big)\; d\gamma :=\Phi_{\mathbf{J}}^\Gamma
\end{equation} with $\mathbf{G}$ the flux for the angular momentum and 
$$
\mathcal{\mathbf{G}}(U^{(p)})=\frac{1}{2}\big ( \mathbf{G}(U^{(p)})+\mathbf{G}(U^{(0)})\big )
$$
or
$$
 \mathcal{\mathbf{G}}(U^{(p)})=\mathbf{G}\bigg ( \frac{U^{(p)}+U^{(0)}}{2}\bigg)
$$
depending on how the system \eqref{eq:1} has been discretized.

Of course, in general, the relation \eqref{cons:J} and \eqref{cons:JJ} cannot be satisfied. In order to be satisfied, as well as keeping 
\begin{equation}
\label{momentum1}\sum_{\sigma \in K} \Phi_{\mathbf{m},\sigma}^K=\oint_K \big( \mathbf{m}^{(p)}-\mathbf{m}^{(0)}\big) d\bx +\Delta t\int_{\partial K} \mathcal{\mathbf{F}}_\mathbf{m}(U^{(p)})\cdot \bn \;d\gamma
\end{equation}
\begin{equation}
\label{momentum2}\sum_{\sigma \in \Gamma} \Phi_{\mathbf{m},\sigma}^\Gamma=\Delta t\int_{\Gamma} \mathcal{\mathbf{F}}_\mathbf{m}(U^{(p)})\cdot \bn \;d\gamma
\end{equation}
with a similar definition of the momentum flux, we will present a perturbation of the momentum residuals. At this level, it is important to provide the quadrature formula. For
\eqref{cons:J} and \eqref{momentum1}, we consider
$$\int_K  f(x) \; d\bx  \approx\frac{|K|}{3}\sum_{\sigma \in K} f_\sigma,$$
so it is exact for \eqref{momentum1} and only approximate for \eqref{cons:J}, but second order.

To achieve this, following \cite{abg}, we introduce a perturbation of the momentum residual, $\br_\sigma^K$, such that the new momentum residual is $\Phi_{\mathbf{m},\sigma}^K+\br_\sigma^K$. We must have:
\begin{equation}
\label{correction}
\begin{split}
\sum_{\sigma\in K}& \br_\sigma^K=0\\
\sum_{\sigma\in K} &\bx_\sigma\wedge \br_\sigma^K = \Phi_{\mathbf{J}}^K-\sum_{\sigma\in K} \bx_\sigma\wedge\Phi_{\mathbf{m},\sigma}^K(U^{(p)}):=\Psi^K
\end{split}
\end{equation}
Since \eqref{cons:JJ} is not true, as above, we introduce  a vectorial correction of the momentum residual $\br_\sigma^\Gamma$. We need:
\begin{equation}
\begin{split}
\sum_{\sigma\in \Gamma}& \br_\sigma^\Gamma=0\\
\sum_{\sigma \in \Gamma} & \bx_\sigma\wedge \br_\sigma^\Gamma =\Phi_{\mathbf{J}}^\Gamma-\sum_{\sigma\in \Gamma} \bx_\sigma\wedge\Phi_{\mathbf{m},\sigma}^\Gamma(U^{(p)}):=\Psi^\Gamma
\end{split}
\end{equation}
$\br_\sigma ^\Gamma$ can be derived in the same way.

In dimension $d$ and for an element with $p$ DOFs, we have $2\times d$ equations for $2\times p$ unknowns. Since $p$ is always larger than $d+1$, the system has a priori solutions. However, finding a close form formula is element dependent. In the following, we show this is possible for triangles and tetrahedrons. Once these solutions are described, we have a scheme that locally preserves the angular momentum $\bJ_\sigma=\bx_\sigma \wedge \bm_\sigma$ and globally conserves
$$\sum_K\frac{|K|}{3}\bigg ( \sum\limits_{\sigma\in K} \bJ_\sigma\bigg )$$
up to boundary terms.

\begin{remark}[Translation invariance]\label{discreteTrans}
The angular momentum is defined after a frame has been defined, and if one makes a translation of vector $\mathbf{a}$, the scheme is translational invariant as in the continuous case.
\end{remark}
\subsection{Solution for triangular elements}
We first have 
$r_1=-r_2-r_3$, so
$$(\bx_2-\bx_1)\wedge r_2+(\bx_3-\bx_1)\wedge r_3=\Psi \in \R$$
We define 
$$r_2=r (\bx_3-\bx_1), \qquad r_3=r(\bx_1-\bx_2), \qquad r\in \R$$ and get
$$
r\bigg ( \det(\bx_2-\bx_1,\bx_3-\bx_1)+\det (\bx_3-\bx_1, \bx_1-\bx_2) \bigg ) =\Psi,
$$
i.e. since $\det(\bx_2-\bx_1,\bx_3-\bx_1)=2|T|$, we have:
\begin{equation}
\label{angular2D}
r=\dfrac{\Psi}{4|T|}, \qquad r_2=r (\bx_3-\bx_1), \qquad r_3=r(\bx_1-\bx_2), \qquad r_1=r (\bx_2-\bx_3).
\end{equation}

\subsection{Solution for tetrahedrons}
Again,
$r_1=-\sum_{j=2}^4 r_j$, so that
$$\sum_{j=2}^4 (\bx_j-\bx_1)\wedge r_j =\Psi.$$
In order to simplify the notations, we introduce $\bx_{ij}=\bx_i-\bx_j$. We are looking for 
\begin{equation*}
\begin{split}
r_2&=\alpha_2\bx_{13}+\beta_2\bx_{14}\\
r_3&=\alpha_3\bx_{12}+\beta_3\bx_{14}\\
r_4&=\alpha_4\bx_{12}+\beta_4\bx_{13}
\end{split}
\end{equation*}
This gives:
\begin{equation*}
\begin{split}
\Psi&=\alpha_2 \bx_{21}\wedge\bx_{13}+\beta_2\bx_{21}\wedge\bx_{14}\\
&\qquad +\alpha_3\bx_{31}\wedge\bx_{12}+\beta_3\bx_{31}\wedge\bx_{14}\\
&\qquad +\alpha_4\bx_{41}\wedge\bx_{12}+\beta_4\bx_{41}\wedge \bx_{13}
\end{split}
\end{equation*}
Then, since $\mathbf{a}\cdot \big (\mathbf{b}\wedge\mathbf{c}\big )=\det(\mathbf{a},\mathbf{b},\mathbf{c})$,
\begin{equation*}
\begin{split}
\Psi\cdot\bx_{21}&= \beta_3   \det(\bx_{21}, \bx_{31}, \bx_{14}) + \beta_4   \det(\bx_{21},\bx_{41},\bx_{13} )\\
\Psi\cdot\bx_{31}&= \beta_2   \det(\bx_{31}, \bx_{21}, \bx_{14}) + \alpha_4 \det(\bx_{31},\bx_{41},\bx_{12} )\\
\Psi\cdot\bx_{41}&= \alpha_2 \det(\bx_{41}, \bx_{21}, \bx_{13}) + \alpha_3 \det(\bx_{41},\bx_{31},\bx_{12} )
\end{split}
\end{equation*}
Now,
\begin{equation*}
\begin{split}
3|T|&= \det(\bx_{21}, \bx_{31}, \bx_{14})=- \det(\bx_{21},\bx_{41},\bx_{13} )=-\det(\bx_{31}, \bx_{21}, \bx_{14})\\
&\qquad =\det(\bx_{31},\bx_{41},\bx_{12} )=\det(\bx_{41}, \bx_{21}, \bx_{13})=-\det(\bx_{41},\bx_{31},\bx_{12} )
\end{split}
\end{equation*}
so that we have:
\begin{equation*}
\begin{split}
\Psi\cdot\bx_{21}&= 3|T|\big (\beta_3  - \beta_4 \big )\\
\Psi\cdot\bx_{31}&= 3|T|\big (-\beta_2   + \alpha_4\big )\\
\Psi\cdot\bx_{41}&= 3|T|\big ( \alpha_2 - \alpha_3 \big ).
\end{split}
\end{equation*}
This suggests to assume $\beta_3 =-\beta_4$, $\alpha_4=-\beta_2$, $\alpha_2=-\alpha_3$, i.e
$$
\begin{array}{c}
\beta_3 =-\beta_4=\dfrac{\Psi\cdot\bx_{21}}{6|T|}\\ ~\\
\alpha_4=-\beta_2=\dfrac{\Psi\cdot\bx_{31}}{6|T|}\\~\\
\alpha_2=-\alpha_3=\dfrac{\Psi\cdot\bx_{41}}{6|T|}
\end{array}
$$
Remembering that $u\wedge(v\wedge w)=(u\cdot w)v-(u\cdot v)w$, we see that
\begin{equation*}
\begin{split}
r_2&=\dfrac{1}{6|T|}\bigg(\Psi\cdot\bx_{41}\;\bx_{13}-\Psi\cdot\bx_{31}\;\bx_{14}\bigg )=\dfrac{1}{6|T|}\Psi\wedge(\bx_{14}\wedge\bx_{31})\\
r_3&=\dfrac{1}{6|T|}\bigg(\Psi\cdot\bx_{21}\;\bx_{14}-\Psi\cdot\bx_{41}\;\bx_{12}\bigg )=\dfrac{1}{6|T|}\Psi\wedge(\bx_{21}\wedge\bx_{14})\\
r_4&=\dfrac{1}{6|T|} \bigg ( \Psi\cdot\bx_{31}\;\bx_{12}-\Psi\cdot\bx_{21}\;\bx_{13}\bigg )=\dfrac{1}{6|T|} \Psi\wedge(\bx_{31}\wedge\bx_{12}).
\end{split}
\end{equation*}
Last, 
\begin{equation*}
\begin{split}
\bx_{14}\wedge\bx_{31}+\bx_{21}\wedge\bx_{14}+\bx_{31}\wedge\bx_{12}&=\bx_{14}\wedge\bx_{31}+\bx_{14}\wedge\bx_{12}+\bx_{31}\wedge\bx_{12}\\
&=\bx_{14}\wedge\bx_{31}+\big ( \bx_{14}+\bx_{31}\big ) \wedge\bx_{12}\\
&=\bx_{14}\wedge\bx_{31}+\bx_{34}\wedge\bx_{12}\\
&=(\bx_{13}+\bx_{34})\wedge\bx_{31}+\bx_{34}\wedge\bx_{12}\\
&=\bx_{34}\wedge\bx_{31}+\bx_{34}\wedge\bx_{12}=\bx_{34}\wedge\bx_{32}
\end{split}
\end{equation*}
and then we have:
\begin{equation}\label{angular:3D}
\begin{split}
r_1&=\dfrac{1}{6|T|}\Psi\wedge\big ( \bx_{23}\wedge\bx_{34}\big )\\
r_2&=\dfrac{1}{6|T|}\Psi\wedge\big (\bx_{14}\wedge\bx_{34}\big)\\
r_3&=\dfrac{1}{6|T|}\Psi\wedge\big (\bx_{14}\wedge\bx_{42}\big)\\
r_4&=\dfrac{1}{6|T|}\Psi\wedge\big ( \bx_{31}\wedge\bx_{32}\big )
\end{split}
\end{equation}

\section{Angular momentum preservation: the high order case}

The idea is similar, the key point is to characterize the quadrature formula that describes $\int_K \bJ\; d\bx$. The additional difficulty is that $f_\sigma=f(\sigma)+O(h^2)$ if we make the geometrical identification of the DOFs with the Greville points.  Here to simplify the notations, $\Phi_\sigma^K$ denotes the residual at $\sigma\in K$ evaluated for the momentum, it does not contain the contribution for the density or the energy. In this section we always use this short hand notation, except at the end.

We start again from \eqref{rd:1} that we rewrite as 
$$U^{(p+1)}_\sigma=U^{(p)}_\sigma+\delta U^{(p)}_\sigma.$$
Then 
$$U^{(p+1)}=\sum_{\sigma} U^{(p+1)}_\sigma \varphi_\sigma.$$
The angular momentum, integrated, is
$$\int_\Omega \bx \wedge \bm^{(p+1)}\; d\bx=\sum_{\sigma} \big (\int_\Omega\bx B_\sigma \; d\bx \big ) \wedge \bm_\sigma^{(p+1)}.$$ We introduce the notations
$$\by_\sigma=\dfrac{1}{|C_\sigma|} \int_\Omega \bx B_\sigma\; d\bx,$$
$$\textbf{z}_\sigma^K=\dfrac{1}{|K|} \int_K \bx B_\sigma\; d\bx,$$
so that we can rewrite the total kinetic momentum as:
$$\int_\Omega \bx \wedge \bm^{(p+1)}\; d\bx=\sum_\sigma |C_\sigma|\; \by_\sigma\wedge \bm_\sigma^{(p+1)},$$ and then we can write the update:
\begin{equation*}
    \begin{split}
 \sum_\sigma |C_\sigma| \;\by_\sigma\wedge \bm_\sigma^{(p+1)}&=
 \sum_\sigma |C_\sigma| \;\by_\sigma\wedge \bm_\sigma^{(p)}- \sum_\sigma \by_\sigma \wedge \big (
 \sum_{K, \sigma\in K} \Phi_{\mathbf{m},\sigma}^K \big )\\
 &=\sum_\sigma |C_\sigma| \;\by_\sigma\wedge \bm_\sigma^{(p)}-\sum_\sigma \by_\sigma \wedge \big (
 \sum_{K, \sigma\in K} \Phi_{\mathbf{m},\sigma}^K +\sum_{\Gamma, \sigma\in \Gamma} \Phi_{\mathbf{m},\sigma}^{\Gamma} \big)\\
  &=\sum_\sigma |C_\sigma| \;\by_\sigma\wedge \bm_\sigma^{(p)}-\sum_\sigma \by_\sigma \wedge \big (
 \sum_{K, \sigma\in K} \Phi_{\mathbf{m},\sigma}^K \big) -\sum_\sigma \by_\sigma \wedge \big (\sum_{\Gamma, \sigma\in \Gamma} \Phi_{\mathbf{m},\sigma}^{\Gamma} \big),
\end{split}
 \end{equation*}
 Since, up to boundary terms, 
 $$\sum_\sigma \by_\sigma \wedge \big (
 \sum_{K, \sigma\in K} \Phi_{\mathbf{m},\sigma}^K \big )=\sum_K \big ( \sum_{\sigma\in K} \by_\sigma\wedge \Phi_{\mathbf{m},\sigma}^K\big ),$$
  $$\sum_\sigma \by_\sigma \wedge \big (
 \sum_{\Gamma, \sigma\in \Gamma} \Phi_{\mathbf{m},\sigma}^\Gamma \big )=\sum_\Gamma \big ( \sum_{\sigma\in \Gamma} \by_\sigma\wedge \Phi_{\mathbf{m},\sigma}^\Gamma\big ),$$
We see that a natural conditions to get local conservation of the kinetic momentum is:
\begin{equation}
\label{angular:HO}
\sum_{\sigma \in K} \by_\sigma \wedge   \Phi_{\bm,\sigma} ^K =\int_K \bx\wedge \big (\bm^{(p)}-\bm^{(0)}\big ) \; d\bx+\Delta t\oint_{\partial K} \mathcal{I}_2\big (\bG(U)\big )\cdot \bn \; d\gamma
\end{equation}
\begin{equation}
\label{angular:H11}
\sum_{\sigma \in \Gamma} \by_\sigma \wedge   \Phi_{\bm,\sigma} ^\Gamma =\Delta t\oint_{\Gamma} \bigg(\hat{\bG}-\mathcal{I}_2\big (\bG(U) \big )\cdot \bn \bigg)\; d\gamma
\end{equation}
If we define the angular momentum as
$$\bJ=\int_\Omega \bx \wedge \bm \; d\bx=\sum_{\sigma}\big ( \int_{\Omega}\bx B_\sigma \; d\bx\big )\wedge \bm_\sigma=\sum_\sigma |C_{\sigma}|\; \by_\sigma \wedge \bm_\sigma,$$ then this quantity is globally conserved.

\bigskip
It is clear that in general, \eqref{angular:HO} is not true, so as in the second order case, we introduce  a vectorial correction of the momentum residual $\br_\sigma^K$, so that modified residuals
$\Psi_\sigma^K=\Phi_{\mathbf{m},\sigma}^K+\br_\sigma^K$ satisfies \eqref{angular:HO}. We need:
\begin{equation*}
\begin{split}
\sum_{\sigma\in K}&\br_\sigma^K=0\\
\sum_{\sigma \in K} &\by_\sigma\wedge \br_\sigma^K =-\sum_{\sigma \in K} \by_\sigma \wedge  \Phi_{\bm,\sigma} ^K \\
&\qquad+\int_K \bx\wedge \big (\bm^{(p)}-\bm^{(0)}\big ) \; d\bx+\Delta t\oint_{\partial K} \mathcal{I}_2\big(\bG(U)\big )\cdot \bn \; d\gamma\\
\end{split}
\end{equation*}
Let us introduce the $\bar\by_K$ barycenter of the $\by_\sigma$ for $\sigma\in K$, 
$$\bar\by_K=\frac{1}{\# \sigma}\sum_{\sigma\in K} \by_\sigma,$$
Here $\# \sigma$ is the number of DOFs in $K$ and we set \footnote{For a 2D vector $\ba=(x,y)$, $\ba^\bot=(-y,x)$, so that $\ba\cdot\ba^\bot=0$ and $\ba\wedge\ba^\bot=||\ba||^2$. in 3D, we have to think a bit.}
$$\br_\sigma^K=\alpha_K  (\by_\sigma -\bar\by_K )^\bot,$$
we get
\begin{equation*}
\begin{split}
\sum_{\sigma\in K}\by_\sigma \wedge \br_\sigma^K=\sum_{\sigma\in K} (\by_\sigma -\bar\by_K )\wedge \br_\sigma^K=
\alpha_K \big [ \sum_{\sigma\in K}  \Vert \by_\sigma-\bar\by_K  \Vert ^2 \big ]=\mathcal{E}
\end{split}
\end{equation*}
with
\begin{subequations}
\label{HO:correction}
\begin{equation}
\label{HO:correction:1}
\mathcal{E}=-\sum_{\sigma \in K} \by_\sigma \wedge  \Phi_{\bm,\sigma}^K +\sum_{\sigma \in K} |K| \textbf{z}_\sigma^K \wedge \big (\bm_\sigma^{(p)}-\bm_\sigma^{(0)}\big ) \; d\bx+\Delta t\oint_{\partial K} \mathcal{I}_2\big(\bG(U)\big )\cdot \bn \; d\gamma
\end{equation}
and
\begin{equation}
\label{HO:correction:2}
\alpha_K \big [ \sum_{\sigma\in K} \Vert \by_\sigma -\bar\by_K \Vert^2 \big ]=\mathcal{E}
\end{equation}
\end{subequations}
Since $\sum\limits_{\sigma\in K} \big \Vert \by_\sigma -\bar\by_K \big \Vert^2\neq 0$ for a simplex, \eqref{HO:correction} always has a solution.

Also in general, \eqref{angular:H11} is not true, we introduce  a vectorial correction of the momentum residual $\br_\sigma^\Gamma$, so that modified residuals
$\Psi_\sigma^\Gamma=\Phi_{\mathbf{m},\sigma}^\Gamma+\br_\sigma^\Gamma$ satisfies \eqref{angular:H11}. We need:
\begin{equation*}
\begin{split}
\sum_{\sigma\in \Gamma}&\br_\sigma^\Gamma=0\\
\sum_{\sigma \in \Gamma} &\by_\sigma\wedge \br_\sigma^\Gamma =-\sum_{\sigma \in \Gamma} \by_\sigma \wedge  \Phi_{\bm,\sigma} ^\Gamma +\Delta t\oint_{\Gamma} \bigg (\hat{\bG}-\mathcal{I}_2\big (\bG(U) \big )\bigg )\cdot \bn \; d\gamma\\
\end{split}
\end{equation*}
Using the same procedure, we can calculate $\br_\sigma ^\Gamma$. Let us note that we need to be consistent with the way that boundary conditions were implemented, so that we preserve conservation at the
boundary for the kinetic momentum. The remark \label{discretTrans} about the translation invariance still applies.

The last thing to do is to define explicitly the vectors $\int_{K}\bx B_\sigma\; d\bx$. We will do it for the polynomial degree k which we consider the cases $k=1,2$ for triangular and quadrilateral elements. First we consider the triangular elements. The case  $k=1$ will give a different solution than that given in the previous section, and it would be interesting to see the difference.

Using the notations of appendix \ref{sec:Bezier}, we write
$$|K|\bx_{k_1k_2k_3}=\int_K \bx B_{k_1k_2k_3} \; d\bx$$ and we note that
$$\bx=\bx_1\lambda_1+\bx_2 \lambda_2+\bx_3 \lambda_3$$
where $\bx_1,\bx_2,\bx_3$ are the vertices of $K$ and $\lambda_i(\bx_j)=\delta_{ij}$. Finally $\lambda_1=B_{100}$, $\lambda_2=B_{010}$, $\lambda_3=B_{001}$ and there is no ambiguity on the degree because the degree of $B_{k_1k_2k_3}$ is $k_1+k_2+k_3$.
\begin{itemize}
\item For $k=1$,
\begin{equation*}
\begin{split}
|K|\;\bx_{100}&=\bx_1\int_K \lambda_1 B_{100}\; d\bx+\bx_2\int_K \lambda_2 B_{100}\; d\bx+\bx_3\int_K  \lambda_3 B_{100}\; d\bx\\
&=\bx_1\int_K  B_{200}\; d\bx+\bx_2\int_K  \frac{1}{2}B_{110}\; d\bx+\bx_3\int_K  \frac{1}{2}B_{101}\; d\bx\\
&=|K|\big (\frac{1}{6}\bx_1+\frac{1}{12} \bx_2+\frac{1}{12} \bx_3\big ) 
\end{split}
\end{equation*}
that is  in the end, for all the DOFs:
\begin{equation*}
\begin{split}
\bx_{100}&=\frac{1}{6}\bx_1+\frac{1}{12} \bx_2+\frac{1}{12} \bx_3 \\
\bx_{010}&=\frac{1}{12}\bx_1+\frac{1}{6}\bx_2+\frac{1}{12}\bx_3\\
\bx_{001}&=\frac{1}{12}\bx_1+\frac{1}{12}\bx_2+\frac{1}{6}\bx_3
\end{split}
\end{equation*}
\item For $k=2$, and the basis functions attached to the 
vertices and the mid-points, we obtain:
\begin{equation*}
\begin{split}
|K|\;\bx_{200}&=\bx_1\int_K \lambda_1 B_{200}\; d\bx+\bx_2\int_K \lambda_2 B_{200}\; d\bx+\bx_3\int_K  \lambda_3 B_{200}\; d\bx\\
&=\bx_1\int_K  B_{300}\; d\bx+\bx_2\int_K  \frac{1}{3}B_{210}\; d\bx+\bx_3\int_K  \frac{1}{3}B_{201}\; d\bx\\
&=|K|\big (\frac{1}{10}\bx_1+\frac{1}{30}\bx_2+\frac{1}{30}\bx_3\big )\\
|K|\;\bx_{110}&=\bx_1\int_K \lambda_1 B_{110}\; d\bx+\bx_2\int_K \lambda_2 B_{110}\; d\bx+\bx_3\int_K  \lambda_3 B_{110}\; d\bx\\
&=\bx_1\int_K  \frac{2}{3}B_{210}\; d\bx+\bx_2\int_K \frac{2}{3} B_{120}\; d\bx+\bx_3\int_K  \frac{1}{3} B_{111}\; d\bx\\
&=|K|\big (\frac{1}{15} \bx_1 + \frac{1}{15}\bx_2+\frac{1}{30}\bx_3\big )\\
\end{split}
\end{equation*}
which gives in the end (for all the DOFs)
\begin{equation*}
\begin{split}
\bx_{200}&=\frac{1}{10}\bx_1+\frac{1}{30}\bx_2+\frac{1}{30}\bx_3\\
\bx_{020}&=\frac{1}{30}\bx_1+\frac{1}{10}\bx_2+\frac{1}{30}\bx_3\\
\bx_{002}&=\frac{1}{30}\bx_1+\frac{1}{30}\bx_2+\frac{1}{10}\bx_3\\
\bx_{110}&=\frac{1}{15} \bx_1 + \frac{1}{15}\bx_2+\frac{1}{30}\bx_3\\
\bx_{101}&=\frac{1}{15} \bx_1 + \frac{1}{30}\bx_2+\frac{1}{15}\bx_3\\
\bx_{011}&=\frac{1}{30} \bx_1 + \frac{1}{15}\bx_2+\frac{1}{15}\bx_3\\
\end{split}
\end{equation*}
\end{itemize}

For quadrilateral elements, we write
$$|K|\bx_{k_1k_2k_3k_4}=\int_K \bx B_{k_1k_2k_3k_4} \; d\bx$$ and 
$$\bx=\bx_1\lambda_1+\bx_2 \lambda_2+\bx_3 \lambda_3+\bx_4 \lambda_4$$
where $\bx_1,\bx_2,\bx_3,\bx_4$ are the vertices of $K$ and $\lambda_i(\bx_j)=\delta_{ij}$. Finally $\lambda_1=(1-x)(1-y)$, $\lambda_2=x(1-y)$, $\lambda_3=xy$, $\lambda_4=(1-x)y$.
For $k=1$, we have
\begin{equation*}
\begin{split}
|K|\;\bx_{1000}&=\bx_1\int_K \lambda_1 B_{1000}\; d\bx+\bx_2\int_K \lambda_2 B_{1000}\; d\bx+\bx_3\int_K  \lambda_3 B_{1000}\; d\bx+\bx_4\int_K  \lambda_4 B_{1000}\; d\bx\\
&=\bx_1\int_K  B_{2000}\; d\bx+\bx_2\int_K \frac{1}{2} B_{1100}\; d\bx+\bx_3\int_K \frac{1}{4} B_{1010}\; d\bx+\bx_4\int_K \frac{1}{2} B_{1001}\; d\bx\\
&=|K|\big (\frac{1}{9}\bx_1+\frac{1}{18} \bx_2+\frac{1}{36} \bx_3+\frac{1}{18} \bx_4\big ) 
\end{split}
\end{equation*}
that is  in the end, for all the DOFs:
\begin{equation*}
\begin{split}
\bx_{1000}&=\frac{1}{9}\bx_1+\frac{1}{18} \bx_2+\frac{1}{36} \bx_3+\frac{1}{18}\bx_4 \\
\bx_{0100}&=\frac{1}{18}\bx_1+\frac{1}{9} \bx_2+\frac{1}{18} \bx_3+\frac{1}{36}\bx_4\\
\bx_{0010}&=\frac{1}{36}\bx_1+\frac{1}{18} \bx_2+\frac{1}{9} \bx_3+\frac{1}{18}\bx_4\\
\bx_{0001}&=\frac{1}{18}\bx_1+\frac{1}{36} \bx_2+\frac{1}{18} \bx_3+\frac{1}{9}\bx_4
\end{split}
\end{equation*}

\section{Test cases}

In this section we will present the numerical results that illustrate the behavior of the residual distribution schemes for the compressible Euler equations that are locally conserving the angular momentum. In the following, we refer to the second order scheme obtained by choosing linear shape functions as B1. Also, Higher order approximation is derived by using quadratic  B\'ezier polynomials (B2) as shape functions. 
\subsection{Isentropic vortex}
\begin{figure}[ht]
\centering
{\includegraphics[width=0.45\textwidth]{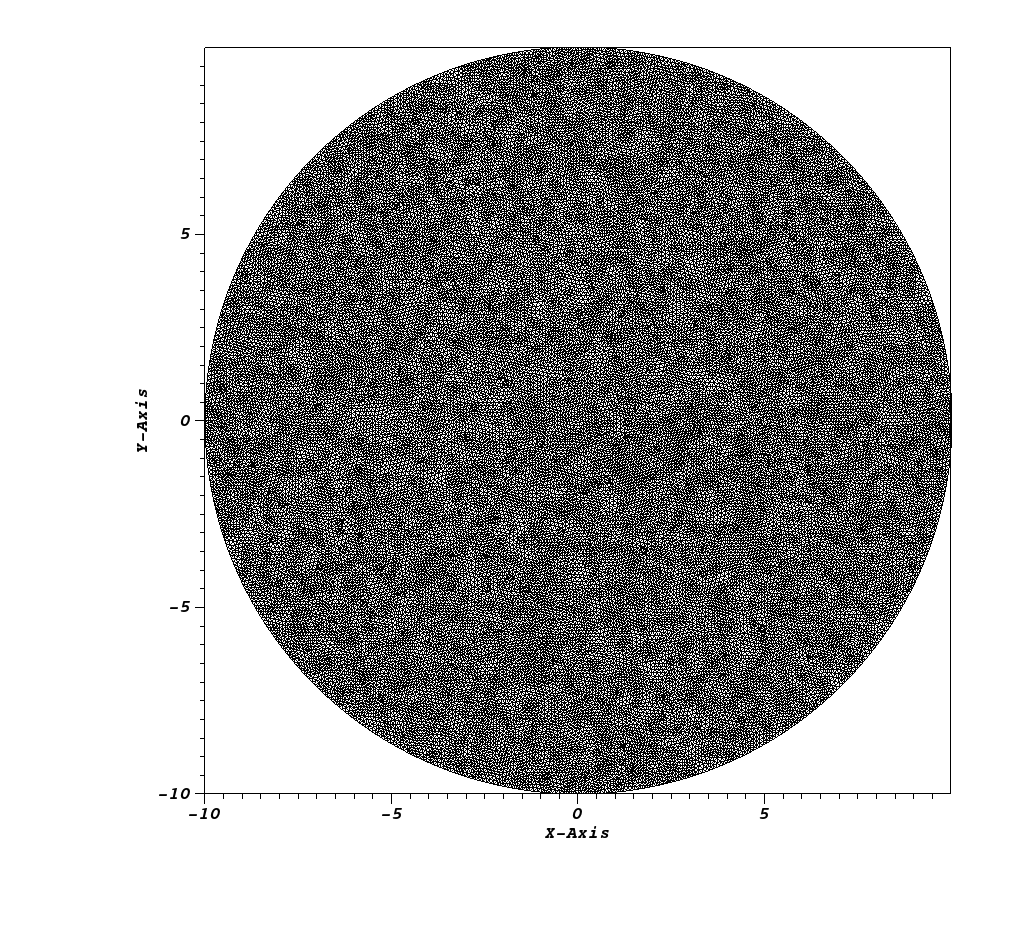}
}
\caption{Mesh for isentropic vortex.}
\label{fig:figure1}
\end{figure}

For the isentropic vortex \cite{yee}, 
we measure the angular momentum throughout the simulation without correction and with correction in the second order and the third order cases on a mesh given by figure \ref{fig:figure1}. 
The physical domain is the circle with radius of 10 and center
at $(x_c,y_c)=(0 , 0)$,  $r=\sqrt{x^2 + y^2}$ and the boundary conditions are \emph{periodic}. The initial conditions for the primitive variables are:
\begin{align*}
\rho &= \left[ 1 - \frac{(\gamma -1)\beta^2}{8\gamma \pi^2}\exp \big( 1 - r^2\big) \right]^{\frac{1}{\gamma-1}}, \\
v_x &= 1 - \frac{\beta}{2\pi}\exp\left(\frac{1-r^2}{2}\right)(y-y_c), \\
v_y &=  \frac{\beta}{2\pi}\exp\left(\frac{1-r^2}{2}\right)(x-x_c), \\
p &= \rho^\gamma,
\end{align*}
for $\gamma = 1.4$, while the free stream conditions are given by:
$$\rho = 1.0, \quad v_{x,\infty} = 1.0, \quad v_{y,\infty} = 0.0, \quad p = 1.0
$$
For all test problems presented in this article,  the reflective wall boundary
conditions are implemented. The final time of the computation is $T=1$. The CFL number is set to $0.5$.

\begin{figure}[ht]
\begin{center}
\subfigure[]{\includegraphics[width=0.45\textwidth]{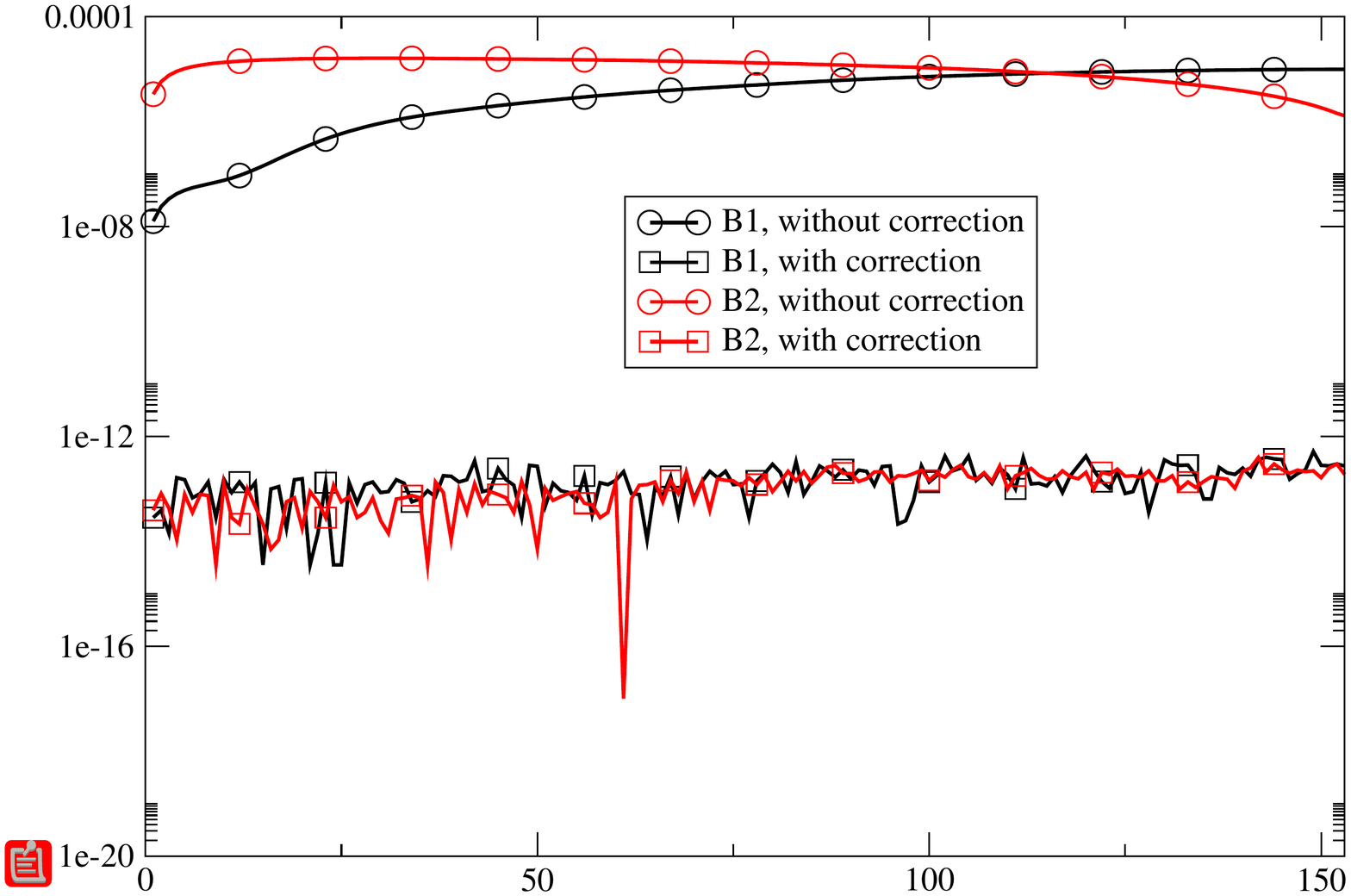}\label{VortexB2vortex-unsteady}}
\subfigure[]{\includegraphics[width=0.45\textwidth]{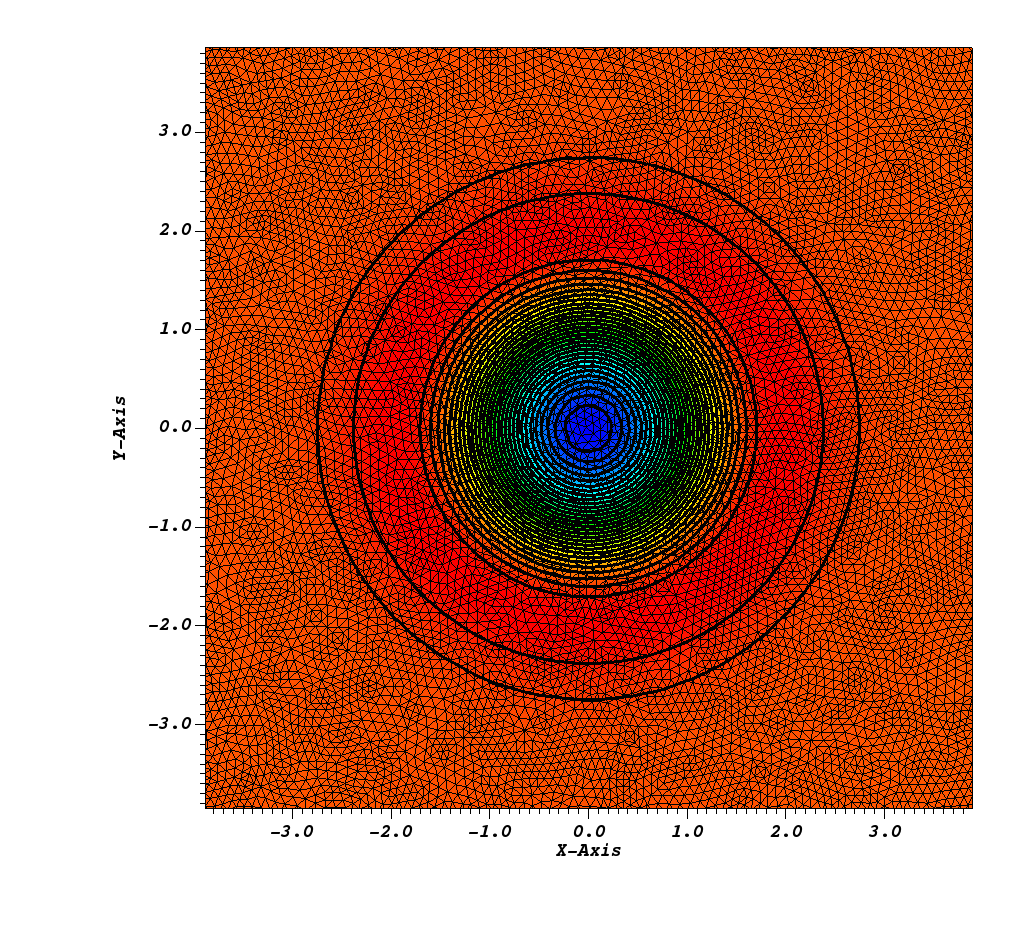}\label{VortexB2vortex-unsteadyy}}
\end{center}
\caption{ (a): Deviation from the initial kinetic momentum for second and third order accurate scheme, with and without correction, (b): solution at time $T=1$, with the mesh}
\end{figure}

In figure \ref{VortexB2vortex-unsteady}, we show the difference between the initial kinetic momentum and the current one for the second order and third order scheme, with and without correction. It is clear that the correction enable to control the kinetic momentum, without negative effect on the solution itself, see figure \ref{VortexB2vortex-unsteadyy} which presents the pressure.

\subsection{Four isentropic vortexes}
\begin{figure}[ht]
\begin{center}
\subfigure[]{\includegraphics[width=0.5\textwidth]{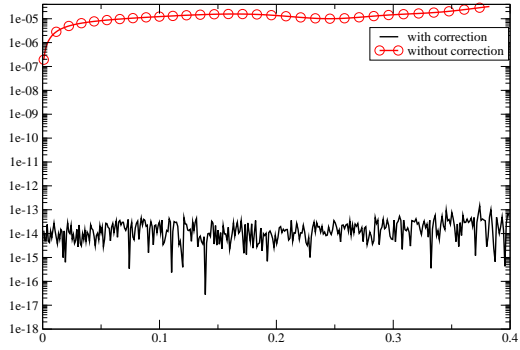}}
\subfigure[]{\includegraphics[width=0.5\textwidth]{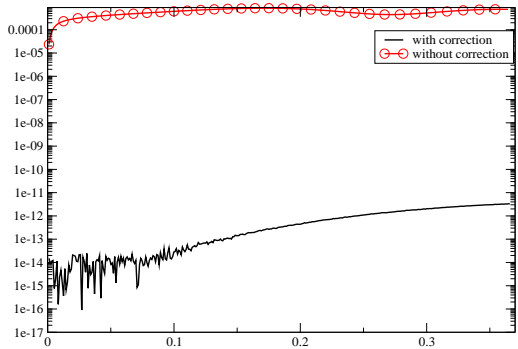}}
\end{center}
\caption{Four isentropic vortices problem, departure  from  the  initial  kinetic  momentum  with and without correction  for: (a) B1, (b) for B2 approximations.}
\label{fig:Scheme4_kineticdifference}
\end{figure}

\begin{figure}[ht]
\begin{center}
\subfigure[B1 with correction]{\includegraphics[width=0.45\textwidth]{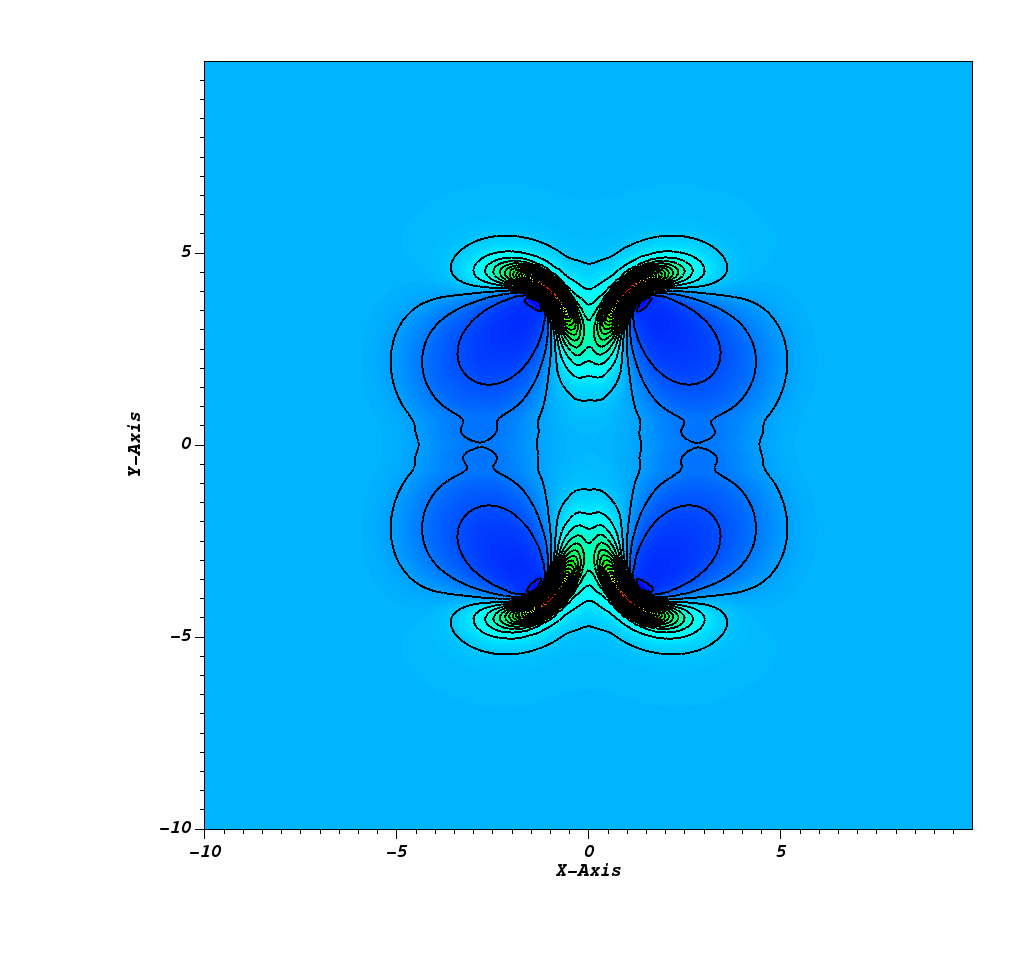}}
\subfigure[B1 without correction.    ]{\includegraphics[width=0.45\textwidth]{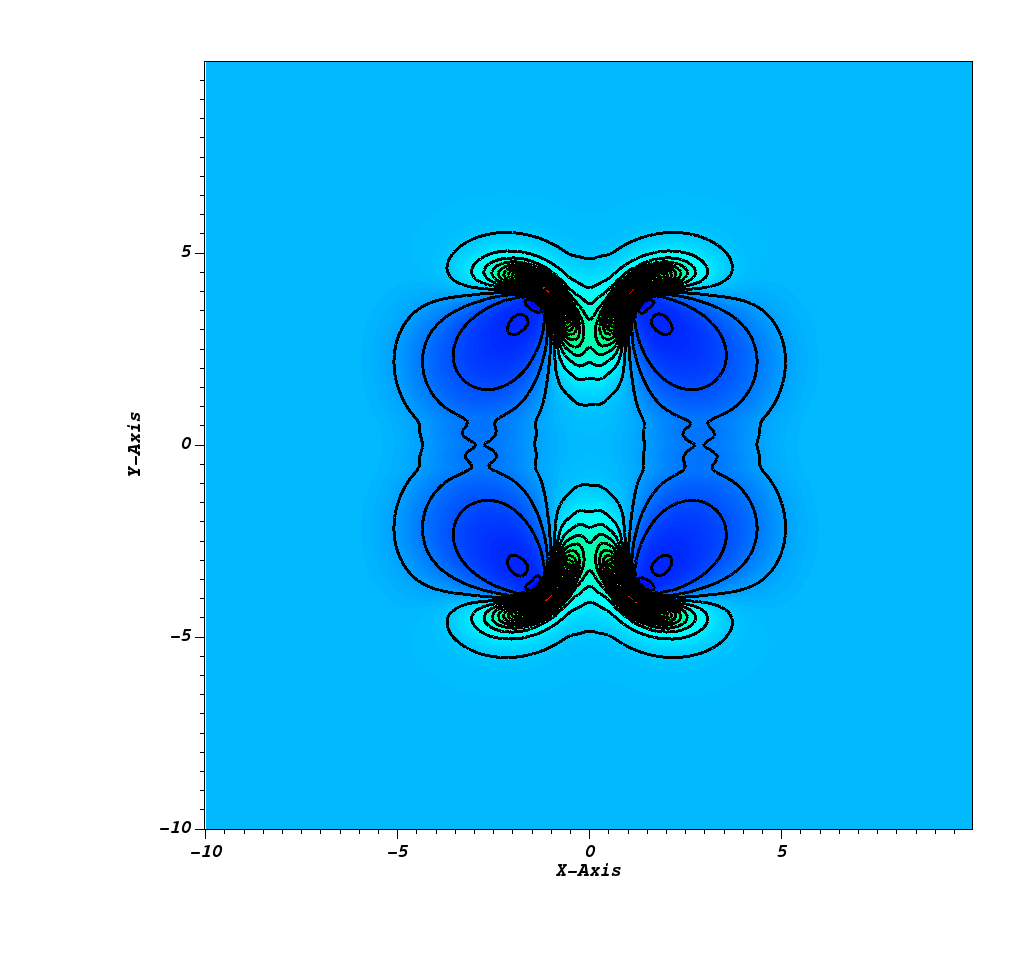}}
\subfigure[B2 with correction]{\includegraphics[width=0.45\textwidth]{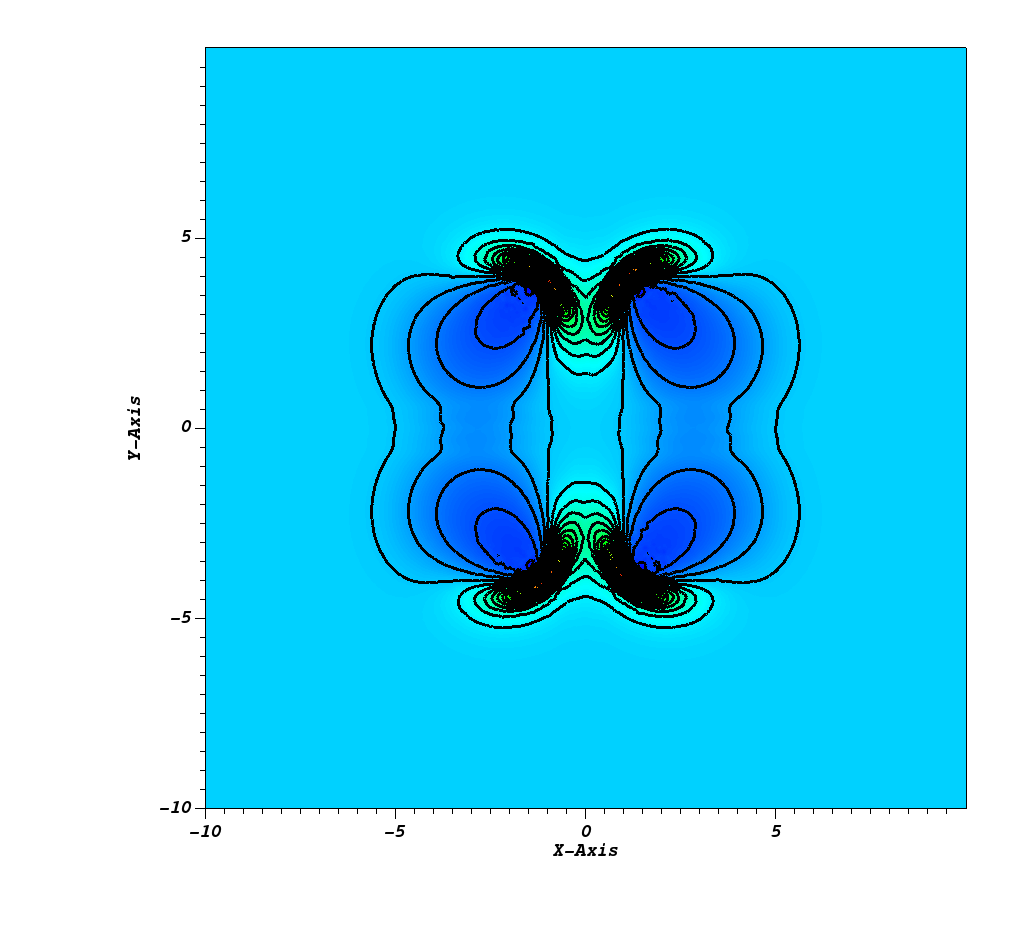}}
\subfigure[B2 without correction.    ]{\includegraphics[width=0.45\textwidth]{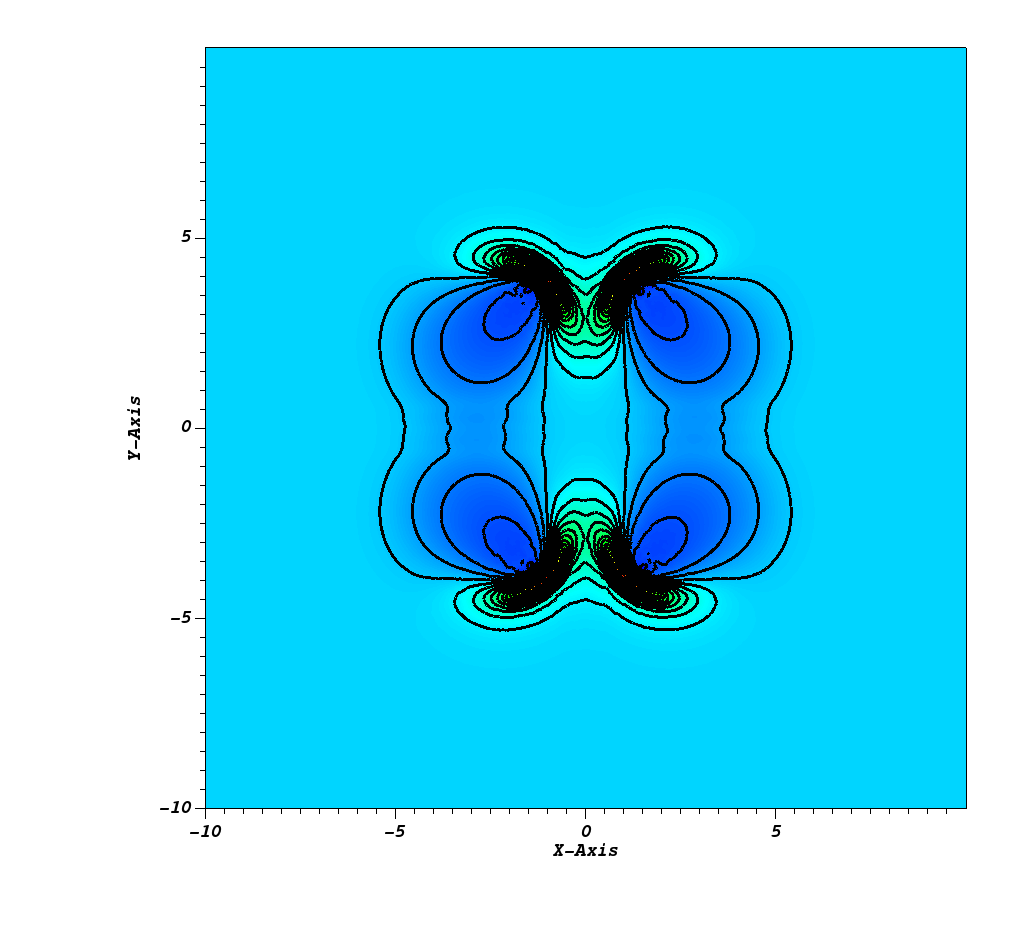}}
\end{center}
\caption{Last solution (pressure) before blow up, $T\approx 0.38$ without correction and $T\approx 0.4$ with correction, using B1 elements, and for B2, $T\approx0.3644$ without correction and $T\approx 0.3657$ with correction.}
\label{fig:Scheme4:p}
\end{figure}

This case is suggested in \cite{gaburro}. 
For this test case, we consider the
four isentropic vortexes centered in $C_1=(2.5,2.5)$, $C_2=(-2.5,2.5)$, $C_3=(-2.5,-2.5)$ and $C_4=(2.5,-2.5)$. The computational domain is a square $[-10, 10]\times[-10,10]$. The initial conditions are given by
\begin{align*}
\rho &= \left[ 1 - \frac{(\gamma -1)\beta^2}{8\gamma \pi^2}\exp \big( 1 - r^2\big) \right]^{\frac{1}{\gamma-1}}, \\
p &= \rho^\gamma,
\end{align*}
\begin{equation*}
\bv =
\left\{
	\begin{array}{ll}
		 \frac{\beta}{2\pi}\exp\left(\frac{1-r^2}{2}\right)(-y,x) & \text{if } xy \geq 0 \\ 
		 \frac{\beta}{2\pi}\exp\left(\frac{1-r^2}{2}\right)(y,-x) & \text{if } xy < 0.
	\end{array}
\right.
\end{equation*}
where
\begin{equation*}
r =
\left\{
	\begin{array}{ll}
		||\bx-C_1||  & \text{if } x \geq 0 \text{ and }y \geq 0 \\
		||\bx-C_2||  & \text{if } x < 0 \text{ and }y \geq 0 \\
		||\bx-C_3||  & \text{if } x < 0 \text{ and }y < 0 \\
		||\bx-C_4||  & \text{if } x \geq 0 \text{ and }y < 0. 
	\end{array}
\right.
\end{equation*}
In figure \ref{fig:Scheme4_kineticdifference} we have represented two different kinetic momentum deviation for B1 and B2 elements: one with regular mesh (obtained from \texttt{gmsh} with the meshing option Frontal Delaunay), and one with the less regular mesh (obtained from \texttt{gmsh} with the meshing option Delaunay). The linear and quadratic meshes have the same number of degrees of freedom. The simulation is obtained using a Galerkin scheme with the CiP stabilisation (stabilisation coefficient set to $0.1$). A priori, this scheme is not adapted because the pressure and the density becomes very small and the gradient of the various variables becomes very large. Indeed, in both cases the scheme  blows up, but the scheme with correction are more robust since the blow up happens later. We have also run the same case with a different numerical strategy that guarantees  positivity preservation  (using a MOOD strategy), and the code does not blow up (we have run until $T=1$).  The pressure field are display on Figure \ref{fig:Scheme4:p}. We note that for $B2$ approximation, the blow up time is slightly larger with correction than without, but the difference is smaller than with B1 approximation.

\subsection{Gresho vortex}
The third considered test case is the Gresho vortex problem, which is a rotating steady solution for the inviscid Euler equations, often used to test conservation of vorticity and angular momentum. 
The angular velocity $v_\phi$ depends only on the radius and the centrifugal force is balanced by the pressure gradient. The physical domain is defined by the circle with radius of 2 and center at $(x_c,y_c) = (0,0)$.
The boundary conditions are gradient free:
\begin{equation*}
\nabla u(\bx)\cdot \bn \rvert_{\bx\in \partial \Omega} = 0, \quad \mbox{for } u \mbox{ a conserved variable } \rho, v_x, v_y, p.
\end{equation*}

The initial conditions for the primitive variables are:
\begin{equation*}
\rho = 1.0, \quad v_x = -v_\phi \frac{(y-y_c)}{r}, \quad v_y = v_\phi \frac{(x-x_c)}{r}, \quad p = p(r),
\end{equation*}
with the orbital velocity $v_\phi$ and pressure $p$:
\begin{equation*}
    v_{\phi}(r) =  \begin{cases}
    5r & r < 0.2 \\
    2-5r & 0.2\leq r < 0.4 \\
    0 & r\geq 0.4
  \end{cases}
\end{equation*}

\begin{equation*}
    p(r) =  \begin{cases}
    5 + \frac{25}{2}r^2 & r < 0.2 \\
    9 - 4\log(0.2) + \frac{25}{2}r^2 - 20r + 4\log(r) & 0.2\leq r < 0.4 \\
    3 + 4\log(2) & r\geq 0.4
  \end{cases}
\end{equation*}

The angular momentum $\vec{J}$ can be written analytically as:
\begin{equation*}
    \vec{J} (r) =  \begin{cases}
    5r^2 & r < 0.2 \\
    2r-5r^2 & 0.2\leq r < 0.4 \\
    0 & r\geq 0.4
  \end{cases}
\end{equation*}
The final time of the computation is $T=0.16$.
and the CFL number is set to 0.25.

\begin{figure}[ht]
\centering
{\includegraphics[width=0.45\textwidth]{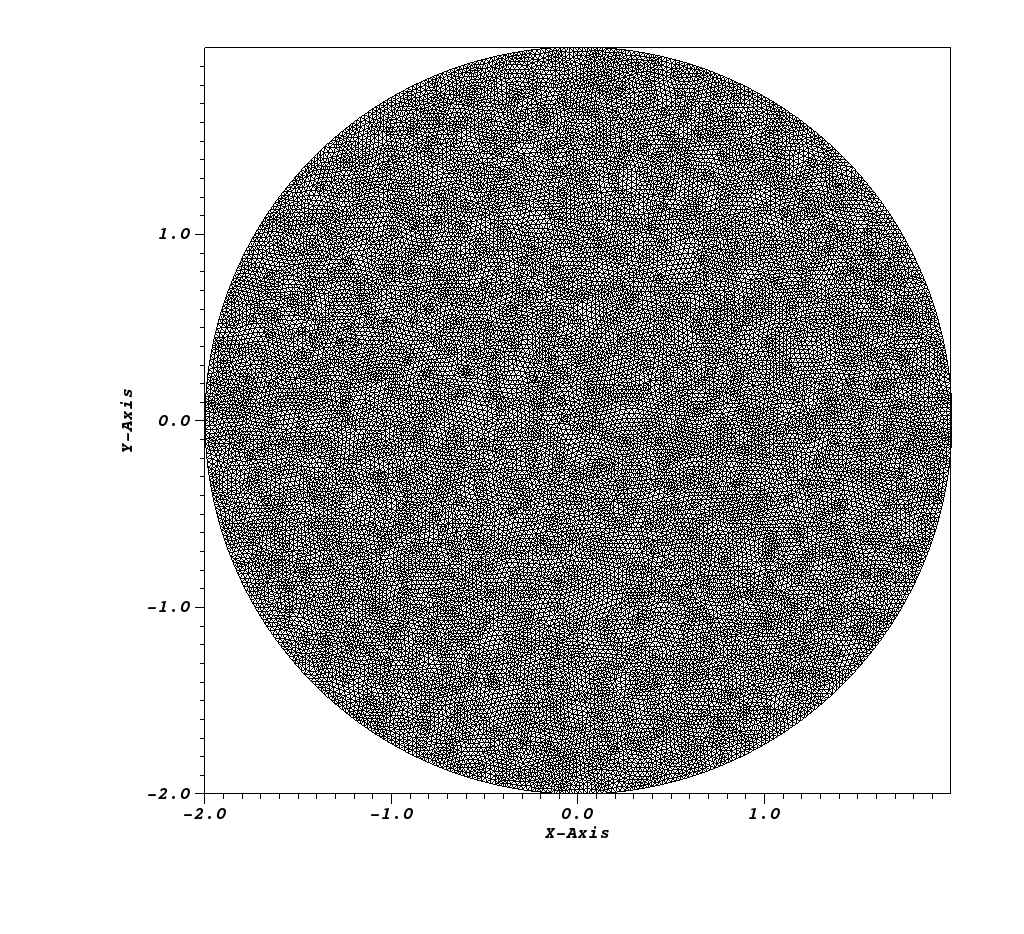}
}
\caption{Mesh for Gresho vortex.}
\label{fig:figure11}
\end{figure}

\begin{figure}[ht]
\begin{center}
\subfigure[PSI]
{\includegraphics[width=0.45\textwidth]{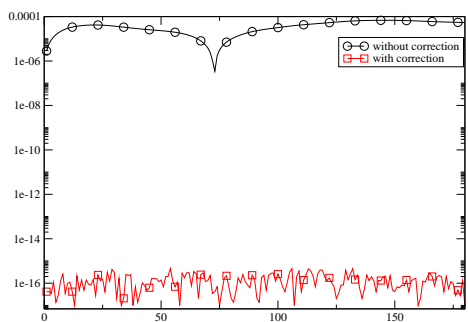}}
\subfigure[Galerkin-CiP]
{\includegraphics[width=0.45\textwidth]{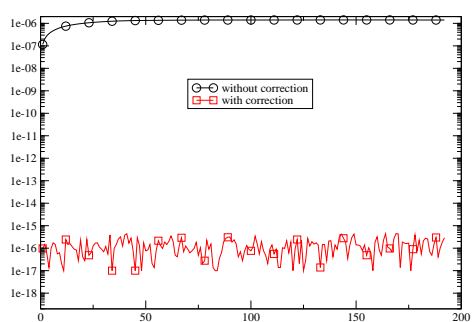}}
\end{center}
\caption{(a): Departure  from  the  initial  kinetic  momentum  with and without correction the PSI scheme with CiP filtering for B2, (b): Departure  from  the  initial  kinetic  momentum  with and without correction the Galerkin scheme with CiP stabilisation for B2.}
\label{fig:figure16}
\end{figure}

\begin{figure}[ht]
\begin{center}
\subfigure[without correction]{\includegraphics[width=0.45\textwidth]{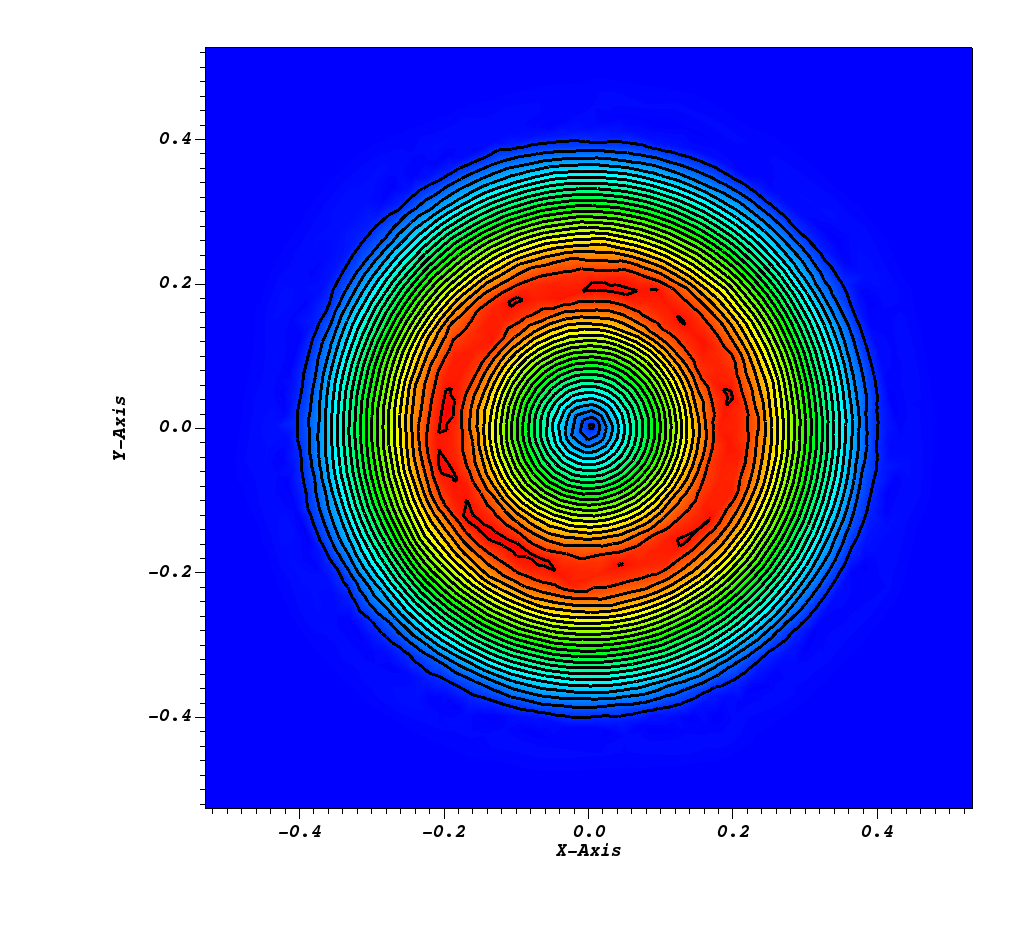}}
\subfigure[with correction.    ]{\includegraphics[width=0.45\textwidth]{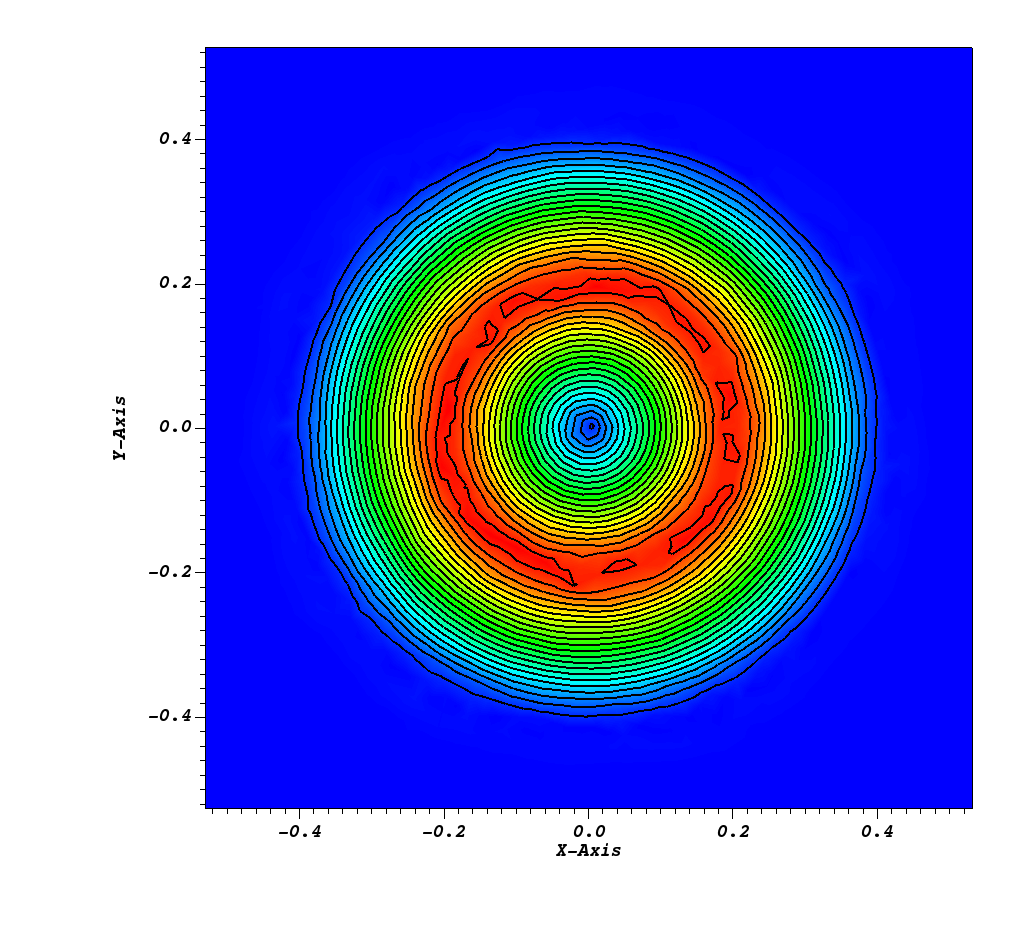}}
\end{center}
\caption{Effect of correction on the velocity field, $T=0.16$ for B2.}
\label{fig:figure17}
\end{figure}

Two schemes are tested, the Galerkin scheme with CiP stabilisation, and the PSI scheme with CiP filtering, see appendix \ref{appendix:PSI} for more details. 
The conservation of kinetic momentum has been tested for the B1 and B2 approximation on a mesh given by figure \ref{fig:figure11}, we only report the results with the B2 approximation, see figure \ref{fig:figure16}, since the results are of similar nature.
The results show clearly the exact conservation of angular momentum obtained with this approach  with B2 elements, the same also hold true for B1 elements.
The obtained solutions at $T=0.16$ almost match in both cases (see figure \ref{fig:figure17}).

\subsection{2D Sod problem}
Further, we have measured the angular momentum on a well-known 2D Sod benchmark problem. The initial conditions
are given by
\begin{equation*}
    (\rho _0,u_0,v_0,p_0) =  \begin{cases}
    ( 1 , 0 , 0 , 1 ),& \text{if } r \leq 0 . 5 , \\
    ( 0 . 125 , 0 , 0 , 0 . 1 ),& \text{otherwise.}
  \end{cases}
\end{equation*}
where $r =\sqrt{ x^2 + y^2}$ is the distance of the point (x , y) from the origin and the computational domain is a square $[-1, 1]\times[-1,1]$.
The final time of the computation is $T=0.16$. 

We have tested two cases. One where the mesh is made of quads, and one when it is made of triangles obtained by cutting the quad into two triangles. The interest of this case is to check if the correction has an influence on the stability property of the initial scheme, and if the structure of the mesh plays an important role. The scheme is formally of the same accuracy as the polynomial approximation, see \cite{larat}, and described in the annex \ref{appendix:PSI}.
The mesh correspond to a $100\times 100$ mesh. 

The pure quad results are displayed in the figure \ref{sod-quad} and those obtained with the triangular mesh are displayed in the next figures.
In figure 
\ref{fig:SodNew_kineticdifference}, we show the evolution of the kinetic momentum over time. In Figure 
\ref{fig:SodNew_rho}, we have the density, with and without correction. In 
figures \ref{fig:SodNew_u} and 
\ref{fig:SodNew_p}, we have displayed the velocity field, and the pressure, with and without correction. We see that the correction has no effect on the stability of the scheme.
\begin{figure}[ht]
\subfigure[]{\includegraphics[width=0.45\textwidth]{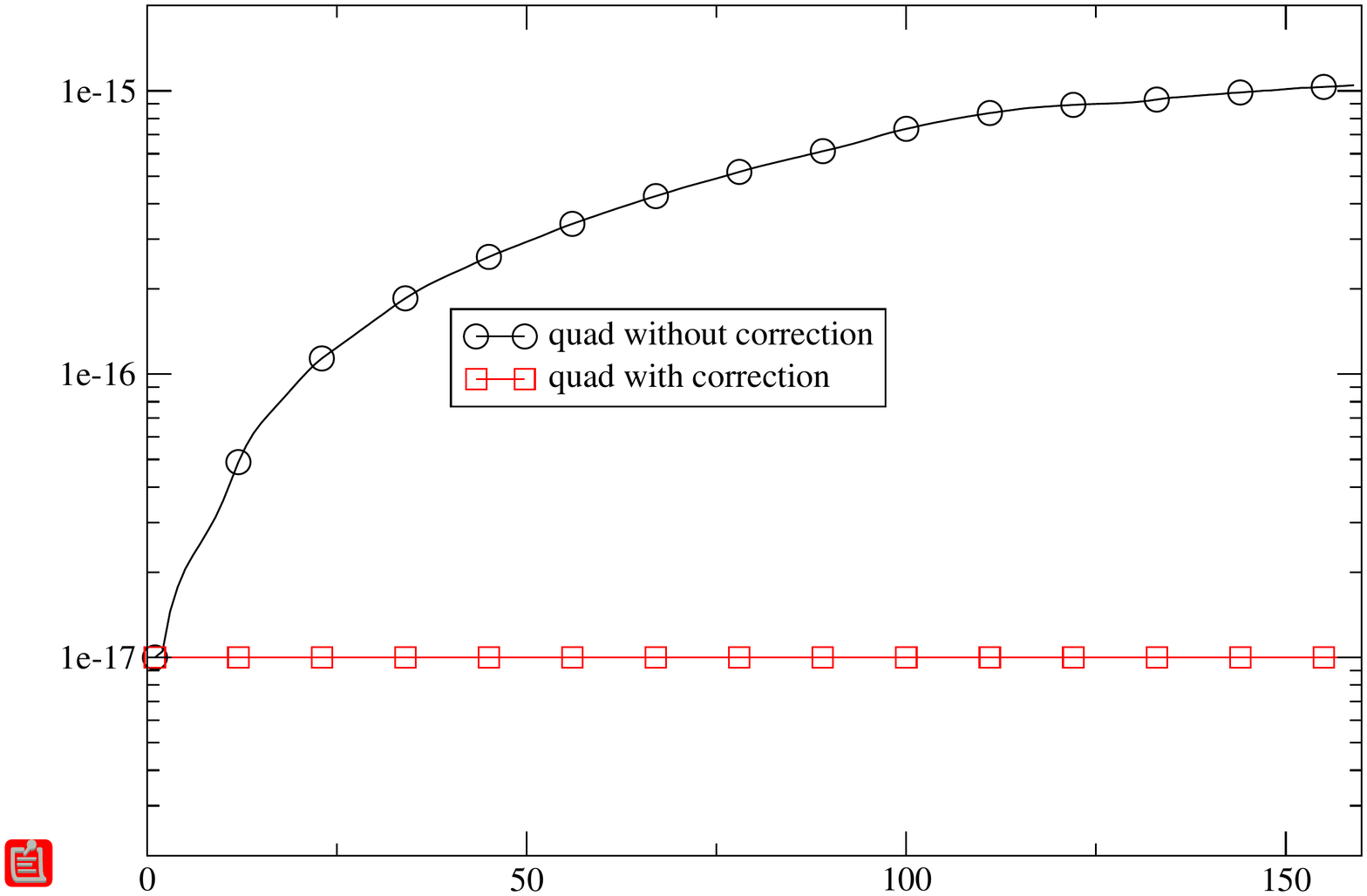}}
\subfigure[]{\includegraphics[width=0.45\textwidth]{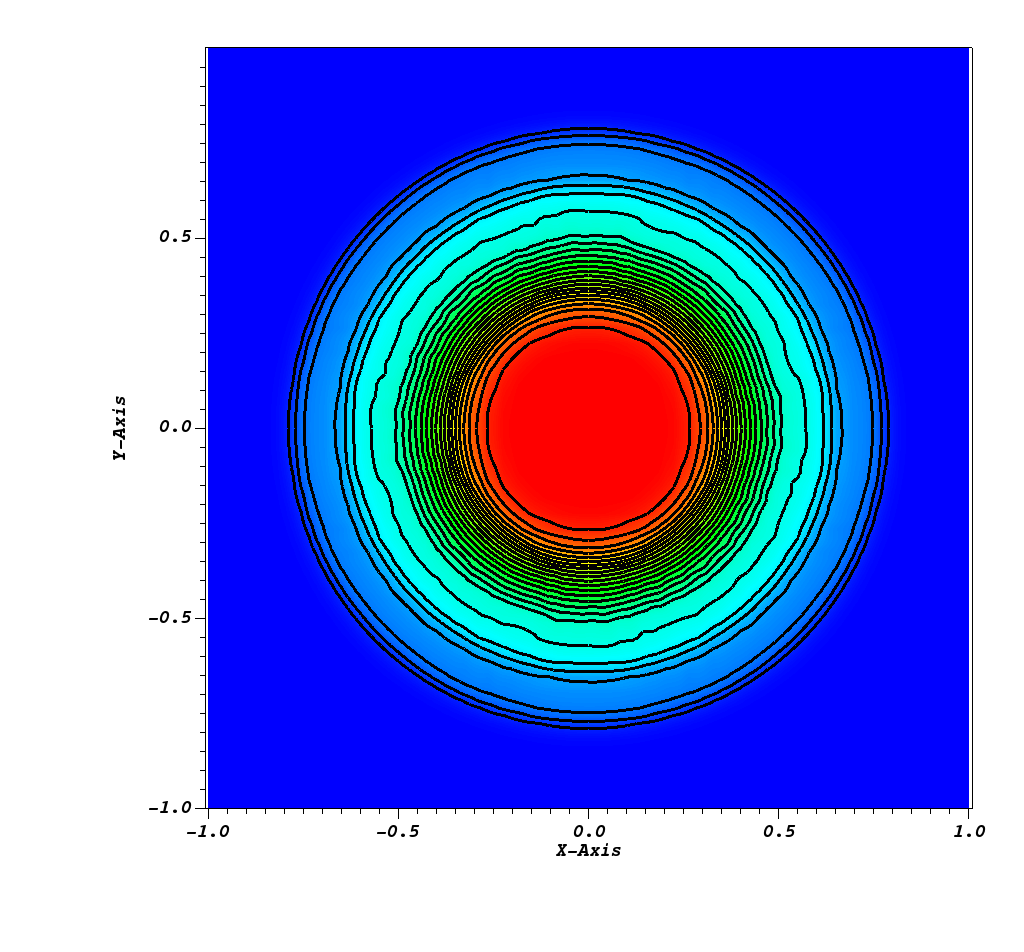}}
\caption{\label{sod-quad} (a): Departure from the initial kinetic momentum with the correction and without the correction, (b): Representation of the density at $T=0.16$, the correction is activated.}
\end{figure}
\begin{figure}[ht]
\begin{center}
{\includegraphics[width=0.45\textwidth]{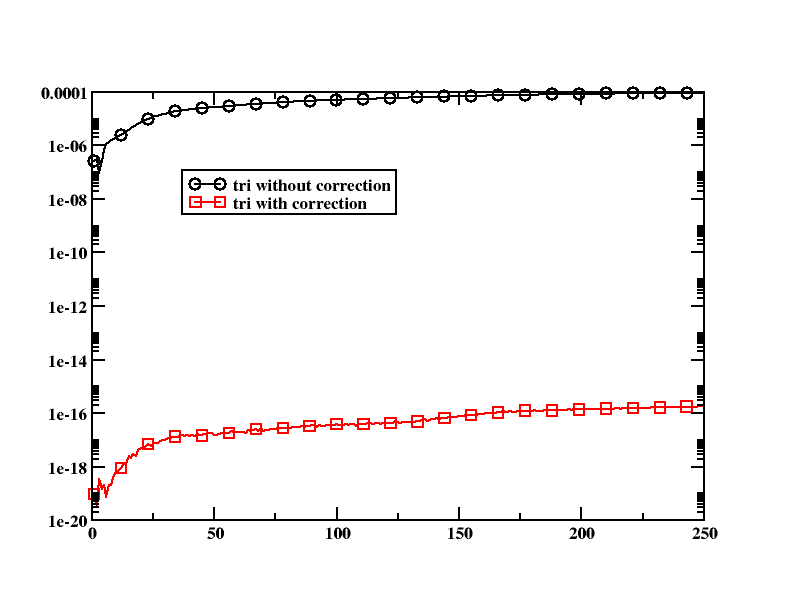}}
\end{center}
\caption{Departure  from  the  initial  kinetic  momentum  with and without correction.}
\label{fig:SodNew_kineticdifference}
\end{figure}

\begin{figure}[ht]
\begin{center}
\subfigure[without correction]{\includegraphics[width=0.45\textwidth]{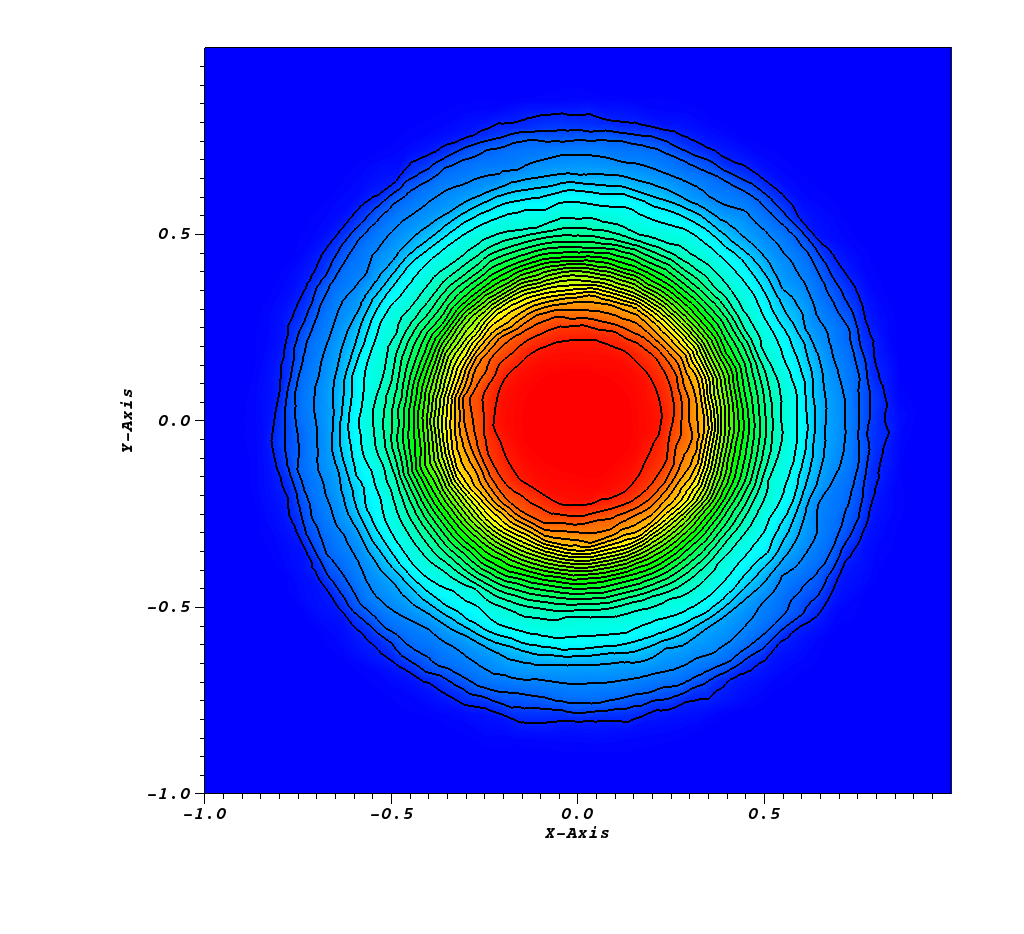}}
\subfigure[with correction.    ]{\includegraphics[width=0.45\textwidth]{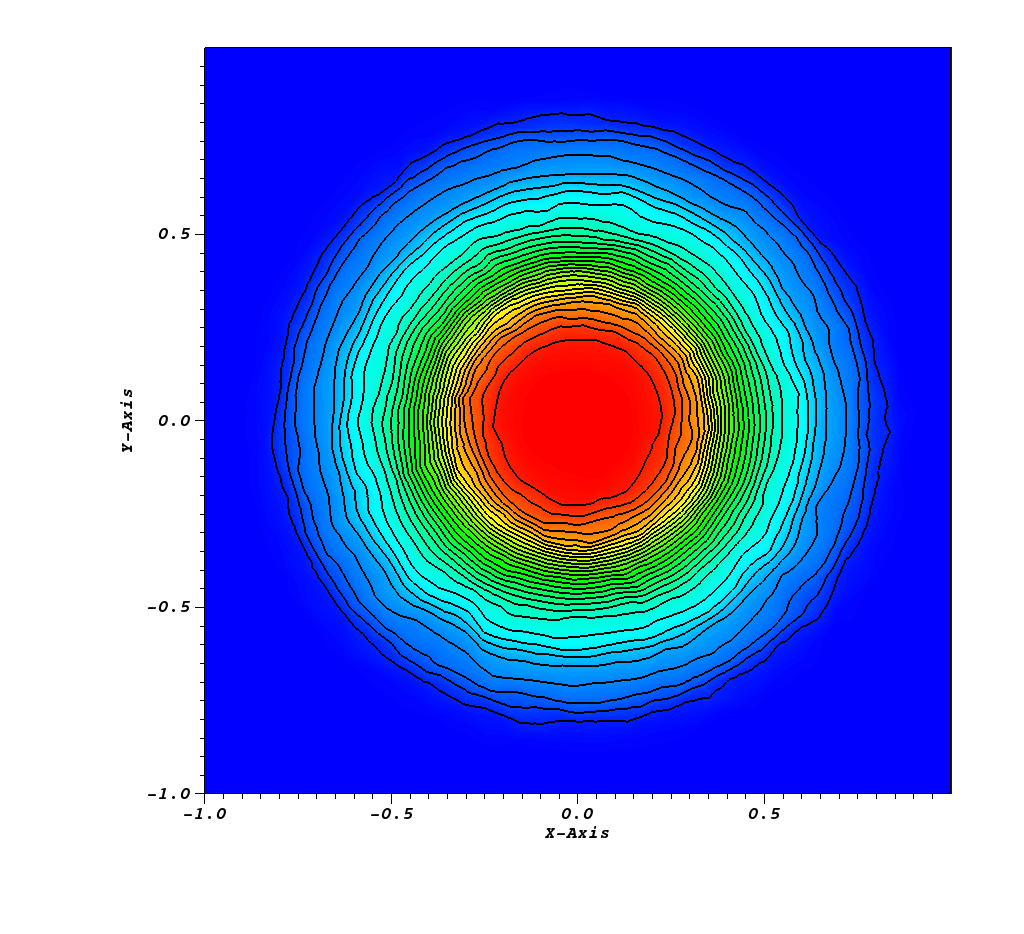}}
\end{center}
\caption{Effect of correction on the density, $T=0.16$ for B2.}
\label{fig:SodNew_rho}
\end{figure}

\begin{figure}[ht]
\begin{center}
\subfigure[without correction]{\includegraphics[width=0.45\textwidth]{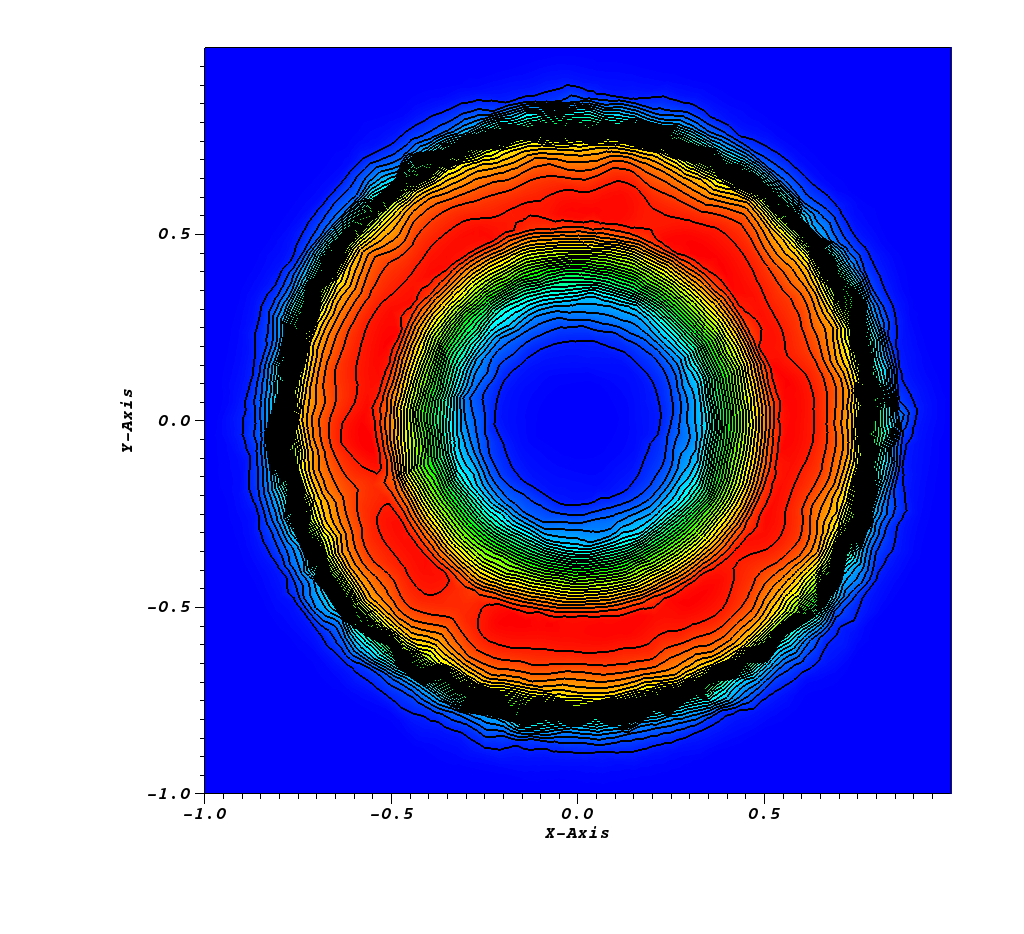}}
\subfigure[with correction.    ]{\includegraphics[width=0.45\textwidth]{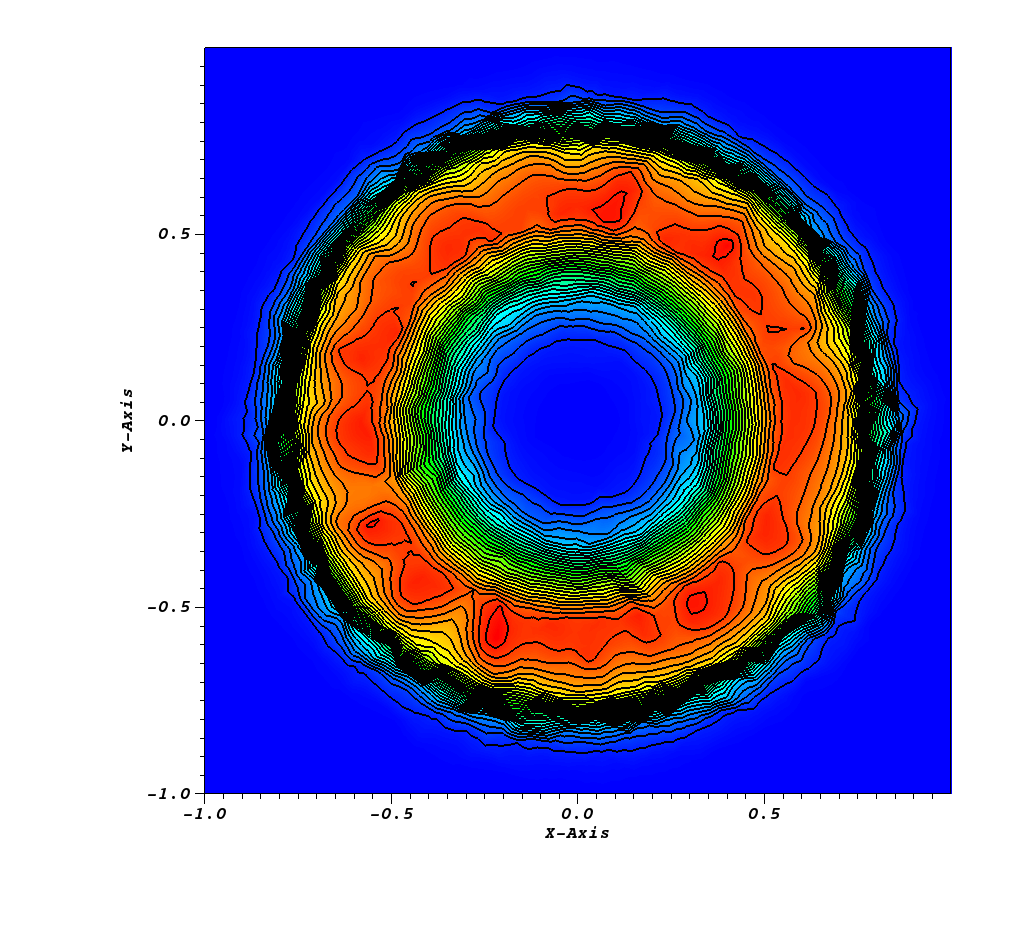}}
\end{center}
\caption{Effect of correction on the velocity field, $T=0.16$ for B2.}
\label{fig:SodNew_u}
\end{figure}

\begin{figure}[ht]
\begin{center}
\subfigure[without correction]{\includegraphics[width=0.45\textwidth]{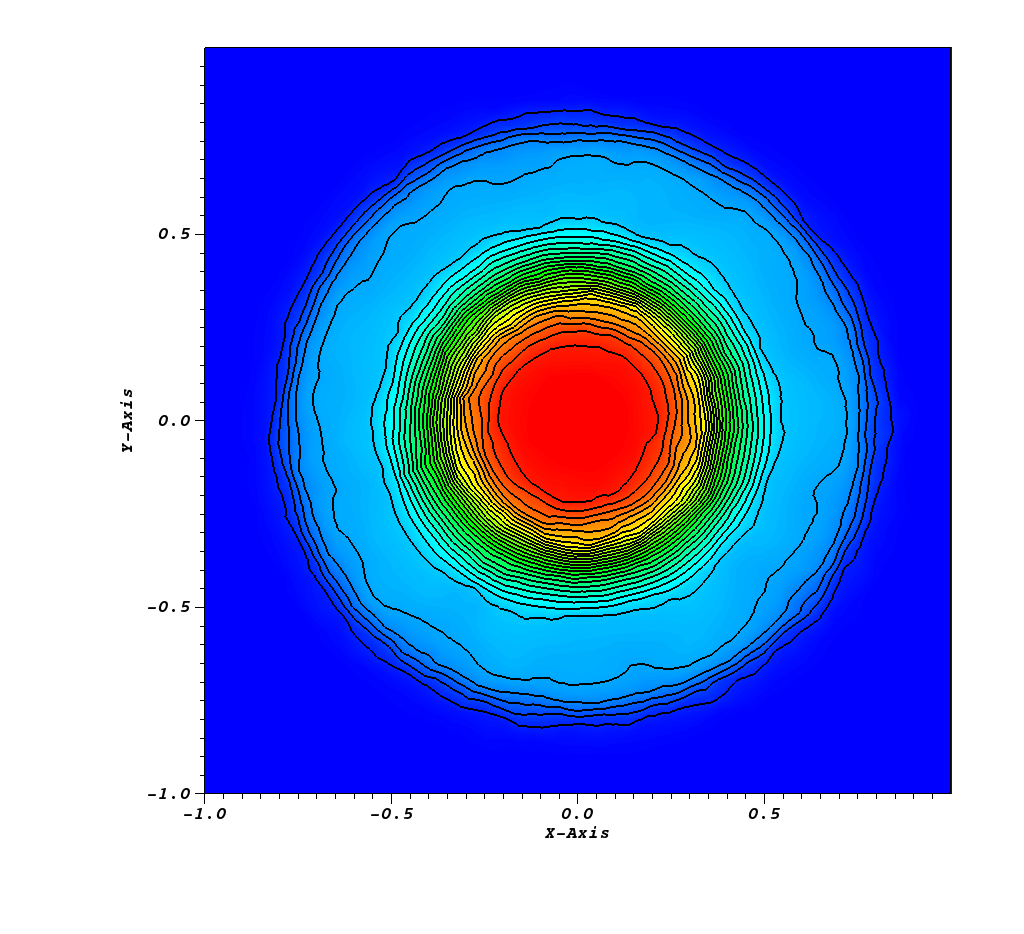}}
\subfigure[with correction.    ]{\includegraphics[width=0.45\textwidth]{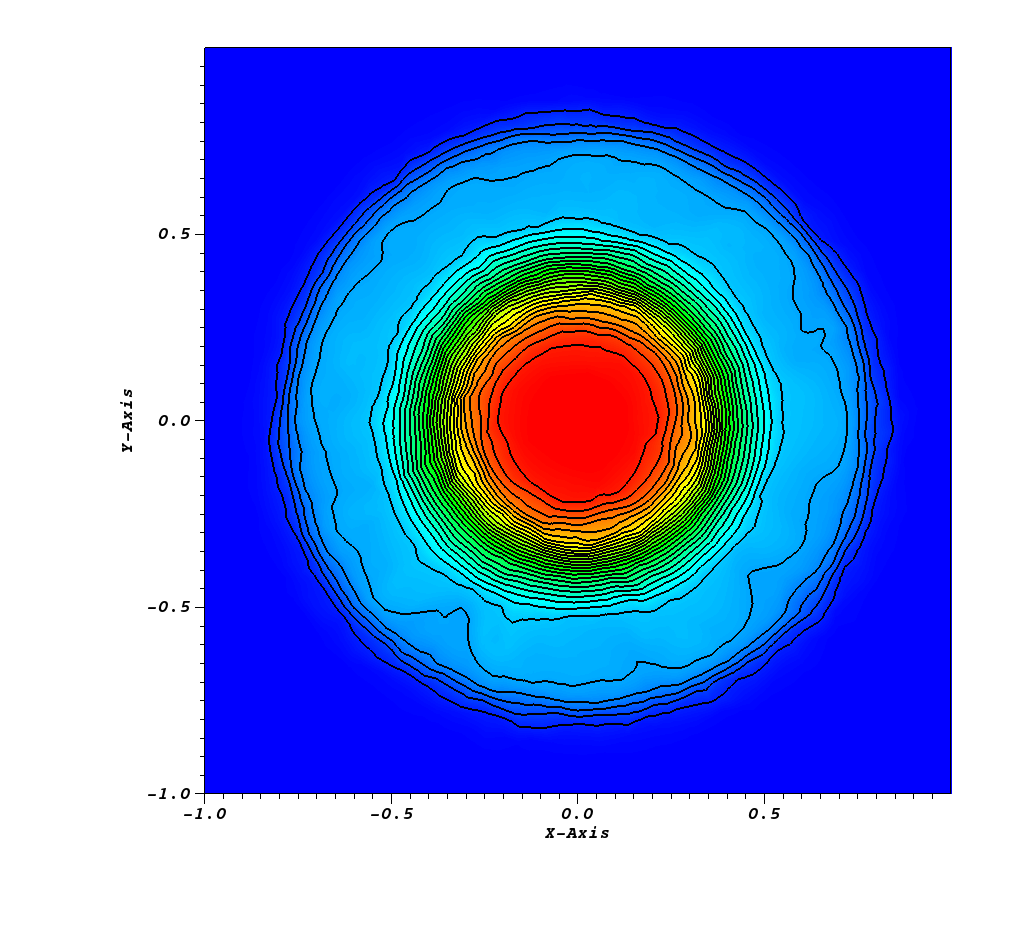}}
\end{center}
\caption{Effect of correction on the pressure, $T=0.16$ for B2.}
\label{fig:SodNew_p}
\end{figure}
We observe that on the pure quad mesh, the correction has a positive effect, though it can also seen  that if it is not active, the variation of the kinetic momentum is negligible. This is in contrast with the triangular mesh, where the effect is much more pronounced. Please note that in both cases, we have the \emph{same} DOFs. It can also be observed that the correction do not have a negative effect on the  non linear stability.
\section{Discussion for DG}

In this section, we discuss how to deal with the problem of kinetic momentum discretization with a DG formulation. The first thing we observe is that the kinetic momentum $\bx \wedge \bm$ and the kinetic momentum flux are obtained from the momentum and the momentum flux simply by multiplying them by polynomials of degree 1 (the space components) and linear combination of these terms. As such, there is nothing special to do, except that we loose systematical one order of accuracy.

The second thing to notice is that the technique developed in this paper can also be applied without any substantial modification, and we keep the same order of accuracy. Let us sketch this.

Using the weak form of \eqref{eq:1}, we get,
\begin{equation*}
M\dpar{\bu}{t}+F=0,    
\end{equation*}
where $M=\big(\int\limits_K \varphi_\sigma \varphi_{\sigma^{'}}\;d\bx\big)_{\sigma,\sigma^{'}\in K}$ is the mass matrix, $\bu=\big( \bu_{\sigma}\big)_{\sigma \in K} ^T$ is the unknown vector and $F=\big(-\int\limits_K \nabla\varphi_{\sigma}\;.\;\bbf(\bu)\;d\bx+\int\limits_{\partial K} \varphi_{\sigma}\bbf(\bu)\;.\;\bn\;d\gamma \big)_{\sigma \in K} ^T$.
Here again,  the  $\varphi_{\sigma}$s are the basis functions.

If we look at the component associated to one DOF, we have
\begin{equation*}
    \sum\limits_{\sigma^{'} \in K} \int\limits_K \varphi_\sigma \varphi_{\sigma^{'}}\dpar{\bu_{\sigma^{'}}}{t}\;d\bx-\int_K \nabla\varphi_{\sigma}\;.\;\bbf(\bu)\;d\bx+\int_{\partial K} \varphi_{\sigma}\bbf(\bu)\;.\;\bn\;d\gamma =0,
\end{equation*}
We can then proceed in the same way as before if the basis functions are B\'ezier polynomials:
\begin{enumerate}
\item Set 
\begin{equation*}
\int_K \varphi_\sigma\frac{\bu_{\sigma}^{(1)}-\bu_{\sigma}^{(0)}}{\Delta t}\;d\bx  -\int_K \nabla\varphi_{\sigma}\;.\;\bbf(\bu^{(0)})\;d\bx+\int_{\partial K} \varphi_{\sigma}\bbf(\bu^{(0)})\;.\;\bn\;d\gamma =0,  
\end{equation*}
\item For any $p\geq 1$, we get
\begin{equation*}
\begin{split}
\int_K \varphi_\sigma\frac{\bu_{\sigma}^{(p+1)}-\bu_{\sigma}^{(p)}}{\Delta t}\;d\bx & + \sum\limits_{\sigma^{'} \in K} \int\limits_K \varphi_\sigma \varphi_{\sigma^{'}} \frac{\bu_{\sigma}^{(p)}-\bu_{\sigma}^{(p-1)}}{\Delta t}\;d\bx -\int_K \nabla\varphi_{\sigma}\;.\; \frac{\bbf(\bu^{(p-1)})+\bbf(\bu^{(p)})}{2}\;d\bx\\
& \qquad  +\int_{\partial K} \varphi_{\sigma}\frac{\bbf(\bu^{(p-1)})+\bbf(\bu^{(p)})}{2}\;.\;\bn\;d\gamma =0.
\end{split}
\end{equation*}
\end{enumerate}
We observe that we have exactly the same formulation as in the globally continuous case, except we have only one residual simply because the set of elements that contain a given degree of freedom is reduced to one element. This being said, we can proceed exactly as in the previous case, without any degradation of the formal accuracy of the method.

\section{Conclusion}
In this paper we have shown, on two example, how to construct systematically  schemes that approximate the compressible Euler equations and are compatible with kinetic momentum preservation. More precisely, starting for  a scheme that is locally conservative, and of formal  of order $r$, one can construct a scheme that is still locally conservative, still formally or order $r$, but also conserves locally the kinetic momentum. The derivation has been done for a residual distribution scheme that assumes a globally continuous approximation of the data, but we have also explain how to extend this to other methods such as discontinuous Galerkin schemes. In the derivation, we have stressed on second order accuracy in time, but using the defect correction approach of \cite{abgrall}, the approach can be easily extended to arbitrary order.

We have illustrated the behavior of the method on several cases, using smooth and non smooth initial conditions.

\section*{Acknowledgments} F.N.M has been funded by the SNF project  	200020\_204917 entitled "Structure preserving and fast methods for hyperbolic systems of conservation laws".
\bibliographystyle{unsrt}
\bibliography{biblio}

\appendix
\section{Examples of fluctuations}\label{appendix:RD}
Here we give two examples which the related residuals satisfy the relevant conservation relations \eqref{conservation:1} or \eqref{conservation:2} \cite{abgrall} depending if we
are considering element residuals or boundary residuals.
\begin{itemize}
    \item the residuals for the SUPG scheme (see \cite{Hughes} for details) are defined by:
\begin{equation*}
\Phi_{\sigma,\bx}^K(\bu)=\int_{\partial K}\varphi_\sigma \bbf(\bbu)\cdot\bn\;d\gamma -\int_K\nabla\varphi_\sigma \cdot \bbf(\bu)\; d\bx+h_K\int_K \bigg(\nabla_{\bu} \bbf(\bbu)\cdot\nabla \varphi_\sigma \bigg)\tau\bigg(\nabla_{\bu}\bbf(\bbu)\cdot\nabla \bu
 \bigg)\; d\bx
\end{equation*}
with $\tau>0$.
\item the residuals for the Galerkin scheme with jump stabilization (see \cite{Burman} for details) are defined by:
\begin{equation*}
\Phi_{\sigma,\bx}^K(\bu)=\int_{\partial K}\varphi_\sigma \bbf(\bbu)\cdot\bn\;d\gamma -\int_K\nabla\varphi_\sigma \cdot \bbf(\bu)\; d\bx+\sum_{e \in K}\theta h_e^2 \int_e [\nabla\bu] \cdot [\nabla \varphi_\sigma
]\; d\gamma
\end{equation*}
with $\theta>0$. Since the mesh is conformal, any internal edge e (or face in 3D) is the intersection of the element K and an other element denoted by $K^{+}$ and for any function $\psi$ we define the jump $[\nabla \psi]=\nabla \psi_{|K}-\nabla \psi_{|K^{+}}$.
\item for the boundary residuals for both cases, we have
\begin{equation*}
\Phi_{\sigma,\bx}^{\Gamma}(\bu)= \int_{\Gamma}   \varphi_\sigma \big( \mathcal{F}_\bn(\bu,g)-\bbf(\bu)\;.\;\bn \big)\;d\gamma.
\end{equation*}
\end{itemize}

\section{PSI scheme}\label{appendix:PSI}
In this appendix, we explain the PSI scheme in more details, in each element $K$, see \cite{Mario,abgrall,paola}. The symbol $\#\sigma$ represents the number of degrees of freedom in the element $K$. First we introduce the Rusanov residuals for a steady version of system \eqref{eq:1}
\begin{equation*}
    \Phi_{\sigma,\bx}^{K,Rus}(\bu)=-\int_K \nabla\varphi_{\sigma}\;.\;\bbf(\bu)\;d\bx+\int_{\partial K} \varphi_{\sigma}\bbf(\bu)\;.\;\bn\;d\gamma+\frac{\alpha_K}{\#\sigma}(\bu_{\sigma}-\bar{\bu}^K),
\end{equation*}
where $\bar{\bu}^K=\frac{1}{\# \sigma}\sum\limits_{\sigma\in K} \bu_\sigma$ and $\alpha_K$ satisfies
$$
\alpha_K\geq  \max_{\sigma,\sigma^{'}\in K}|\int_K \varphi_{\sigma}\nabla \varphi_{\sigma^{'}}\;.\;\nabla_{\bu}\bbf(\bu)\;d\bx|.
$$
Hence, we can write the residual in p-th iteration for \eqref{eq:1} as
\begin{equation*}
  \Phi_{\sigma}^K(U^{(p)})=\Phi_{\sigma}^{K,Rus}(U^{(p)})+\frac{|K|}{\# \sigma} (U_{\sigma}^{(p)}-U_{\sigma}^{(0)}), 
\end{equation*}
We then consider the quasi-linear form of the \eqref{eq:1} in two dimensions
\begin{equation*}
    \dpar{\bu}{t}+A(\bar{\bu})\dpar{\bu}{x}+B(\bar{\bu})\dpar{\bu}{y}=0,
\end{equation*}
where $A(\bar{\bu})$ and $B(\bar{\bu})$ are the Jacobians of the fluxes evaluated at some average state $\bar{\bu}$. Let us introduce a direction $\bd=(\bd_x,\bd_y)$. In this work, we
have chosen to use $\bd=\frac{\bv}{||\bv||}$ since this approach makes the scheme rotationally invariant, and is probably more natural. We also have considered the matrix $T_{\bd}=A(\bar{\bu})\;\bd_x+B(\bar{\bu})\;\bd_y$, which is diagonalizable. So, we can consider $L_i, i=1,4$ and $R_i, i=1,4$ its left and right eigenvectors, respectively. Now we define the following fluctuations by projecting the first order nodal residuals onto a space of left eigenvectors
$$
\Psi_{\sigma}^i=L_i\;.\;\Phi_{\sigma}^K,
$$
We obviously have
$$
\sum_{\sigma \in K} \Psi_{\sigma}^i=L\;.\;\Phi^K:=\Psi^L.
$$
In order to obtain the high order nodal limited residuals, we would compute the distribution coefficients $\beta_{\sigma}^i$ as
$$
\beta_{\sigma}^i=\frac{\max(\frac{\Psi_{\sigma}^i}{\Psi},0)}{\sum\limits_{j=1}^4 \max(\frac{\Psi_{\sigma}^j}{\Psi^L},0)},
$$
and we note that if $\Psi\neq 0$, then $\sum\limits_{j=1}^4 \max(\frac{\Psi_{\sigma}^j}{\Psi^L},0)\geq 1$ so that there is no problem of division.
To this end, the high order nodal residuals are projected back to the physical space
\begin{equation*}
 \Phi_{\sigma}^{*}=\sum\limits_{i=1}^4 (\beta_{\sigma}^i\;\Psi^L)R_i .  
\end{equation*}

\section{B\'ezier polynomials: notations}\label{sec:Bezier}
If $K$ is a simplex, we will denote its vertices as $v$, or $v_i$, or $i$, knowing that we have $d+1$ vertices. The barycentric coordinates with respect to the vertices will be denoted by $\lambda_v$ or $\lambda_i$, depending on the context. The B\'ezier polynomials of degree $n$ are labelled according to a  multi-index with $d+1$ components, $(k_1, \ldots , k_{d+1})$ with $\sum_{i=1}^{d+1} k_i=n$, or a DOF $\sigma$ according to the context. To fix ideas, let us detail the 2D case. The B\'ezier polynomial of index $(k_1,k_2,k_3)$  corresponds to the DOF $\sigma$ that we identify to the point $\bx_\sigma$ in $K$ which barycentric coordinates are $(\tfrac{k_1}{n}, \tfrac{k_2}{n}, \tfrac{k_3}{n})$ which are called the Greville points. We have
$$B_{k_1,k_2,k_3}:=\dfrac{k_1! k_2!k_3!}{(k_1+k_2+k_3)!}\lambda_1^{k_1}\lambda_2^{k_2}\lambda_3^{k_3}:=B_\sigma$$
and we see that
\begin{itemize}
\item 
$\sum_{\sigma\in K}B_\sigma=\sum_{k_i\geq 0, k_1+k_2+k_3=n} B_{(k_1,k_2,k_3)}=1$,
\item $B_\sigma\geq 0$ on $K$,
\item
$$\int_K B_{k_1,k_2,k_3}(\bx)\; d\bx=\dfrac{2}{(n+1)(n+2)} |K|.$$
\end{itemize}
The B\'ezier polynomials of degree $n$ constitute a basis of $\P^n$, the set of polynomials of degree less or equal to $n$. 

We have:
\begin{itemize}
\item $k=1$, $B_{100}=\lambda_1$, $B_{010}=\lambda_2$, $B_{001}=\lambda_3$,
\item $k=2$, $B_{200}=\lambda_1^2$, $B_{020}=\lambda_2^2$, $B_{002}=\lambda_3^2$, $B_{110}=2\lambda_1\lambda_2$,  $B_{101}=2\lambda_1\lambda_3$, $B_{011}=2\lambda_2\lambda_3$, 
\item $k=3$, 
$B_{300}=\lambda_1^3$, $B_{030}=\lambda_2^3$, $B_{003}=\lambda_3^3$, $B_{210}=3\lambda_1^2\lambda_2$, $B_{120}=3\lambda_1\lambda_2^2$, $B_{021}=3\lambda_2^2\lambda_3$, $B_{012}=3\lambda_1\lambda_2^2$, $B_{201}=3\lambda_1^2\lambda_3$, $B_{102}=3\lambda_1\lambda_3$, $B_{111}=6\lambda_1\lambda_2\lambda_3$.
\end{itemize}
\end{document}